\documentclass[review]{elsarticle}

\biboptions{sort&compress}
\usepackage{epsfig, graphicx}
\usepackage{latexsym,amsfonts,amsbsy,amssymb}
\usepackage{amsmath}
\usepackage{cleveref}
\usepackage{multirow}
\usepackage{tikz}
\usepackage{subcaption}
\usepackage{mathtools}
\usepackage{booktabs}
\usetikzlibrary{shapes,arrows,positioning,decorations.pathreplacing,
        intersections,matrix,patterns,fadings,calc,decorations.markings,arrows.meta}
\usepackage{pgfplots}
\usepackage{url}

\newcommand{\vect}[1]{\textbf{#1}}
\newcommand{\vu}{\vec{\vect{u}}}
\newcommand{\vf}{\vec{\vect{f}}\,}
\newcommand{\vk}{\vec{k}}
\newcommand{\vbk}{\vec{\textbf{k}}}
\newcommand{\vp}{\vec{p}\,}
\newcommand{\dx}[0]{\, \mathrm{d} \mathbf{x}}

\usepackage{todonotes}

\begin{document}

\begin{frontmatter}

\title{Monolithic multigrid for implicit Runge--Kutta discretizations of incompressible fluid flow\tnoteref{t1}}
\tnotetext[t1]{Submitted to the editors DATE.}
\tnotetext[t2]{This research was enabled in part by support provided by ACENET (www.ace-net.ca),
Scinet (www.scinethpc.ca), and Compute Canada (www.computecanada.ca). The work of SM and RA
was partially supported by an NSERC Discovery Grant. The work of PEF was supported by the UK Engineering and Physical Sciences Research Council [grant numbers EP/R029423/1 and EP/W026163/1].}

\author[1]{Razan Abu-Labdeh\corref{cor1}}
\ead{rabulabdeh@mun.ca}

\author[1]{Scott MacLachlan}
\ead{smaclachlan@mun.ca}

\author[2]{Patrick E.~Farrell}
\ead{patrick.farrell@maths.ox.ac.uk}

\cortext[cor1]{Corresponding author}
\address[1]{Department of Mathematics and Statistics, Memorial University of Newfoundland, St.~John's, NL, Canada}
\address[2]{Mathematical Institute, University of Oxford, Oxford, UK}

%\author{Razan Abu-Labdeh\thanks{Department of Mathematics and Statistics, Memorial University of Newfoundland, St.~John's, NL, Canada 
%	(\email{rabulabdeh@mun.ca}, \email{smaclachlan@mun.ca}).}
%  \and Patrick E.~Farrell\thanks{Mathematical Institute, University of Oxford, Oxford, UK
%		(\email{patrick.farrell@maths.ox.ac.uk}).}
%  \and Scott MacLachlan\footnotemark[2]}

\begin{abstract}
Most research on preconditioners for time-dependent PDEs has focused on
implicit multi-step or diagonally-implicit multi-stage temporal discretizations. In this paper, we consider
monolithic multigrid preconditioners for fully-implicit multi-stage
Runge--Kutta (RK) time integration methods. These temporal discretizations have
very attractive accuracy and stability properties, but they couple the spatial
degrees of freedom across multiple time levels, requiring the solution of
very large linear systems. We extend the classical Vanka
relaxation scheme to implicit RK discretizations of saddle point problems.
We present numerical results for the incompressible Stokes, Navier--Stokes,
and resistive magnetohydrodynamics equations, in two and three
dimensions, confirming that these relaxation schemes lead to robust
and scalable monolithic multigrid methods for a challenging range of
incompressible fluid-flow models.
\end{abstract}


\end{frontmatter}

%\linenumbers

\section{Introduction}

Among the many applications of advanced computer simulation, models of fluid flow have been a persistent and common driving force in research and practice.  The history of spatial discretization of fluid problems dates back at least to the 1960's (e.g., the MAC-scheme discretization of the Navier--Stokes equations \cite{FHHarlow_JEWelch_1965a, AJChorin_1968a,RTemam_1969a}), but continues to this day with investigation of higher-order mixed finite-element discretizations for both Newtonian and complex fluids \cite{scott1985conforming, SZhang_2005a, VJohn_etal_2017a, NRGauger_etal_2019a, DSchotzau_2004a, KHu_etal_2017a, KHu_JXu_2019a, KHu_etal_2020a}.
Alongside this thrust to higher-order spatial discretizations comes a need for stable higher-order temporal discretizations, for which implicit Runge-Kutta methods are a natural choice. % However, this exposes a weakness in the current literature, where effective linear and nonlinear solution algorithms for these temporal discretizations are somewhat critically under-studied.
%For many of these problems, coupling the spatial discretization with implicit time-stepping is attractive, due to its improved stability over explicit timestepping, but this requires the development of effective nonlinear and linear solution algorithms for the resulting time-stepping equations.  Furthermore, while much is known about higher-order temporal discretizations for both ODEs and PDEs, many implicit numerical methods for time-dependent fluid simulations still make use of low-order temporal discretization, such as implicit Euler or the second-order backward differentiation formula (BDF2).
In this paper, we investigate the development of efficient Newton--Krylov--multigrid strategies for implicit Runge--Kutta discretizations of incompressible fluid-flow problems.

Effective solver strategies for both stationary problems and time-dependent flow models discretized via either multi-step schemes or diagonally implicit Runge-Kutta (DIRK) schemes have been studied for many years.  For time-dependent Newtonian flows, both fully and semi-implicit pressure-correction schemes (e.g., \cite{AJChorin_1968a, AJChorin_1969a, RTemam_1969a, Kan86, JBBell_PColella_HMGlaz_1989a}) have been proposed, based primarily on multigrid solution of the pressure-Poisson equation, but the construction and analysis of general high-order schemes is non-trivial \cite{JLGuermond_etal_2006a}.  Monolithic multigrid schemes (both linear and nonlinear) have also been broadly considered, first arising in the late 1970's and early 1980's \cite{ABrandt_NDinar_1979a, ABrandt_1984a}.  More approaches have been proposed since these early works, including techniques for Newtonian flows based on Vanka \cite{vanka1986block} and Braess--Sarazin \cite{braess1997efficient} relaxation, and generalizations of these techniques to more complex flow settings and discretizations \cite{DGstokes, adler2016monolithic, scottmhd}.  Simultaneously, block preconditioning strategies have also been developed, for a variety of discretizations and flow settings \cite{DKay_etal_2002a, elman, Wathen_etal_2017a, PEFarrell_etal_2019a, PEFarrell_etal_2021a, laakmann2021augmented}.  Despite this substantial body of work on multi-step methods, there are (to our knowledge) few comparable publications on solution strategies for multi-stage implicit Runge--Kutta (IRK) discretizations of flow models~\cite{pazner2017stage, southworth2021fast1}. 

A small body of work exists on solvers for IRK discretizations for parabolic PDEs~\cite{van2005multigrid, ERosseel_etal_2008a, TBoonen_etal_2009a, mardal2007order, LJay_2000a, HChen_2014a, rana2020new, southworth2021fast1, southworth2021fast2}.  Much of this work focuses on block-structured preconditioners for the tensor-product systems generated by IRK discretization \cite{mardal2007order, LJay_2000a, HChen_2014a, rana2020new} where, for example, standard multigrid methods can be used to solve the diagonal blocks.  The recent method of Southworth et al.~\cite{southworth2021fast1, southworth2021fast2} appears to be very effective, again leveraging standard preconditioners for linear systems corresponding to BDF-type discretizations.  On the other hand, the work of Vandewalle and others \cite{van2005multigrid, ERosseel_etal_2008a, TBoonen_etal_2009a} applies monolithic multigrid methods to these discretizations, using block-Gauss--Seidel type relaxation for parabolic equations and a block-extension of the Hiptmair relaxation \cite{RHiptmair_1999a} for the eddy-current form of the curl-curl equation.  Similar block-Jacobi relaxation was used for both the heat and Gross--Pitaevskii equations in \cite{Irksome}.  Here, we investigate extensions of Vanka relaxation for IRK discretizations of fluid flow problems.

In this paper, we consider standard mixed finite-element (spatial) discretizations of Stokes, Navier--Stokes, and magnetohydrodynamic (MHD) flows, coupled with IRK discretizations in time.  We focus on the development of monolithic geometric multigrid preconditioners for the coupled systems of equations to be solved at each timestep.  For nonlinear problems, we use these preconditioners in a standard Newton--Krylov--multigrid setting, using Newton's method to linearize the coupled nonlinear systems at each timestep.  We expect the same techniques would apply to the various simplifications of Newton's method that are applicable in the IRK context \cite{butcher1976implementation, bickart1977efficient}. Numerical results are presented for standard benchmarks in two and three spatial dimensions, showing that this solution approach is equally effective for IRK discretizations as it is for BDF discretizations. 

The remainder of this paper is organized as follows.  In \Cref{sec:RKdisc}, we review the Runge--Kutta discretization approach for systems of ODEs.  For fluid-flow models, this is typically used in a method-of-lines approach with some spatial discretization, and \Cref{sec:spatial_disc} reviews mixed finite-element discretization of the Stokes, Navier--Stokes, and MHD models considered here.  In \Cref{sec:monolithicMG}, we present the constituent parts of the monolithic multigrid algorithm that we propose for solution of the resulting linear(ized) systems of equations.  Numerical results that confirm the effectiveness of this approach are given in \Cref{sec:numerics}.  Finally, conclusions and directions for future work are given in \Cref{sec:conclusion}.
  
\section{Runge--Kutta temporal discretizations}
\label{sec:RKdisc}
While BDF (and other linear multi-step) schemes can achieve
higher-order convergence, they do so at a cost to their stability, with
the widely known result that no linear multi-step scheme with order
greater than two can be A-stable (the so-called Second Dahlquist
Barrier) \cite{wanner1996solving}.  Because of this (and other reasons),
Runge--Kutta integrators are widely used when we seek higher-order time
integration methods.  In contrast to multi-step schemes (where solutions at past time-steps are used in the approximation),
Runge--Kutta methods are \textit{multi-stage} schemes, where a number
of intermediate stage values are used to achieve the approximation. In general, an $r$-stage Runge--Kutta method applied to the system of ordinary differential equations $u'(t)=f(u(t),t)$ is given by 
 \begin{equation}
\begin{aligned}
k_{i}&=f\left( u^{n}+\Delta t\sum_{j=1}^{r}a_{ij}k_{j} ,t^{n}+c_{i}\Delta t\right),\text{ for }i=1,2,\ldots,r,\\
u^{n+1}&=u^{n}+\Delta t \sum_{j=1}^{r} b_{j}k_{j}.
\end{aligned}
\label{general RK}
\end{equation}
The coefficients in the scheme are the stage
times (or nodes) $c_i$, the weights $b_j$, and the Runge--Kutta
matrix $A = [a_{ij}]$.  Taken together, these form the Butcher tableau for a given scheme~\cite{butcher1976implementation, JCButcher_2006}. For consistency, we require that
$\sum_{j=1}^{r}b_{j}=1$ and $\sum_{j=1}^{r} a_{ij}=c_{i}$, for all
$i=1 ,2, \dots ,r$. The $r$ stage values are represented by the set
$\{k_{i}\}_{i=1}^r$ and the approximation at time $t^{n} = t^0 + n\Delta t$ is denoted by $u^{n}$.

Runge--Kutta methods are generally classified by the non-zero pattern of the matrix $A$.  Methods can be explicit, with $a_{ij}=0 \text{ } \forall j \geq i$, or implicit, when $\exists j \geq i$ with $a_{ij}\neq 0$. The implicit methods can further be classified into diagonally implicit, with $a_{ij}=0 \text{ } \forall j > i$, or fully implicit, when $\exists j > i$ such that $a_{ij} \neq 0$.  Further specialization is also possible, such as singly diagonally implicit Runge--Kutta (SDIRK) methods, which are diagonally implicit (DIRK) methods with the added property that $a_{ii} = a_{jj}$ for all $i$ and $j$, and explicit singly diagonally implicit (ESDIRK) methods, which have an all-zero first row of $A$, followed by SDIRK structure on lower rows (of which the Crank--Nicolson scheme is a well-known example). 

There are three main points to consider when choosing a Runge--Kutta method, regarding its stability, accuracy, and computational cost.
For any scheme, we define the function $r(z)$ as the map produced when applying the scheme to the (scalar, linear) Dahlquist test problem, $u' = \lambda u$ for $\lambda\in\mathbb{C}$, with $u^{n+1} = r(\lambda\Delta t)u^{n}$.  The domain of stability of the scheme is defined as the region in the complex plane where $|r(z)|<1$.  RK methods are said to be \textit{A-stable} if their domain of stability includes the entire left-half of the complex plane.  If, additionally, we have that $\lim_{z\rightarrow -\infty} |r(z)| = 0$, we say that the scheme is \textit{L-stable}.  For many applications, L-stability is the preferred property, since an L-stable scheme generally damps non-physical high-frequency oscillations that may pollute a numerical solution.  As is typical, explicit Runge--Kutta (ERK) methods have finite regions of stability, and only implicit Runge--Kutta (IRK) schemes can be A- or L-stable.

The local truncation error of an RK scheme is defined as the error made in a single step of the scheme, starting with the analytical solution of the differential equation as $u^n$, compared to $u(t^{n+1})$, while the global error is the accumulated error in the approximate solution over the timesteps needed to reach a fixed time.  We typically discuss such errors by their order, meaning that we bound the error by a constant (depending on $f(u,t)$ and the analytical solution, $u(t)$) times $(\Delta t)^p$ to establish that a scheme has order $p$.  Typically (e.g., when $f(u,t)$ is continuous in $t$ and Lipschitz continuous in $u$), the global error is one order less than the local truncation error.  A well-known result is that the order of global error of an ERK method cannot be greater than its number of stages (and, to achieve order $p>5$, an ERK scheme must have at least $p+1$ stages) \cite[Section 324]{butcher2016numerical}.  In contrast, the maximum order of global error for an IRK discretization can be as much as twice the number of stages in the scheme.  While higher-order global error is attractive, for both stiff DEs and systems of differential-algebraic equations (DAEs), the so-called \textit{stage order} of a Runge--Kutta method is more important \cite[Section 362]{butcher2016numerical}.  Here, in addition to the truncation error, the accuracy of a scheme is determined by also bounding the approximation of stage $i$ to $u(t^{n} + c_{i}\Delta t)$ by some constant (depending on $f(u,t)$ and $u(t)$) times $(\Delta t)^{q+1}$, thus defining the stage order as $\text{min}\{q,p\}$. For index-2 DAEs (as are considered here), the order of accuracy of a scheme is limited by its stage order, due to perturbation bounds on the solution of the constrained system~\cite[Section VII.4]{wanner1996solving}.  This greatly limits our choice of schemes that allow higher-order accuracy.  While DIRK methods can have reasonable global order, their stage order is typically limited to 1. We note that ESDIRK methods are an exception to this, with stage order limited to 2, due to the structure of their Butcher tableau. In contrast, the stage order of fully IRK schemes can be as large as the number of stages, making these the preferred schemes for integrating DAEs.

The downside of IRK schemes is their computational cost.  ERK methods can be implemented at the cost of one evaluation of $f(u,t)$ for each stage in the method.  In contrast, IRK methods require solution of a system of equations for each timestep (that may be large when $u$ represents a spatially discretized approximation to the solution of a PDE).  Herein lies the attraction of DIRK, SDIRK, and ESDIRK schemes.  In these approaches, rather than having to solve for the stages in a coupled manner, each stage can be solved for sequentially, allowing the reuse of standard linear and nonlinear solvers from backward-Euler type schemes.  SDIRK and ESDIRK afford even more of an advantage, particularly in the linear case, as the same solver architecture can be directly reused in the solution process for each stage.  General IRK methods, in contrast, do not allow this simplification.  While block-preconditioning strategies can be used to again leverage existing solver architectures from the multistep case \cite{mardal2007order, LJay_2000a, HChen_2014a, rana2020new, southworth2021fast1, southworth2021fast2}, these theoretical results tend to be limited to simple cases, excluding (for example) nonlinear systems of DAEs, as arise in standard models of computational fluid dynamics.

In this paper, we consider a standard Newton--Krylov--multigrid framework for the solution of the nonlinear systems of equations that arise from using general IRK discretizations for the Navier--Stokes equations and the equations of magnetohydrodynamics.  Because the details of these solvers depend directly on the spatial discretization, we next discuss the mixed finite-element discretization of these models.

\section{Discretization of fluid models}
\label{sec:spatial_disc}

In this section, we consider the interplay of mixed finite-element spatial discretization for incompressible models of fluid flow with temporal discretization by IRK methods.  We consider three models: the linear Stokes model, the (nonlinear) Navier--Stokes equations,  and the equations of single-fluid visco-resistive incompressible magnetohydrodynamics (MHD).  In~\Cref{sec:monolithicMG}, we will focus on the development of a monolithic multigrid methodology for the linearized systems that result from applying Newton's method to the nonlinear problems.  Both here and in that exposition, we will focus on the details of the algorithm for the simplest case of the linear Stokes model.

\subsection{Time-dependent Stokes equations}
In the viscous limit of incompressible flow, inertial forces in the model can be neglected, leading to the time-steady Stokes equations.  We consider here the time-dependent analogue of the Stokes equations
on a bounded Lipschitz domain $\Omega \subset \mathbb{R}^d$, $d \in \{2, 3\}$:
\begin{subequations}
\label{time stokes}
\begin{alignat}{3}
\rho\textbf{u}_{t}-\mu \Delta \textbf{u}+\nabla p&=\textbf{f}     && \text{  in } \Omega \times (0,T_{f})\label{subeqn:velocity}\\
-\nabla \cdot \textbf{u}&=0     && \text{  in } \Omega \times (0,T_{f}),\\
\textbf{u}&=0     &&\text{  on } \partial\Omega \times (0,T_{f}),\\
\textbf{u}(\textbf{x},0)&=\textbf{g}(\textbf{x})      &&\text{  on } \Omega \times \{ t=0 \},
\end{alignat}
\end{subequations}
where $\textbf{u}(\textbf{x},t)$ is the velocity, $p(\textbf{x},t)$ is the pressure, and $\textbf{f}(\textbf{x},t)$ is a suitably smooth forcing term. Here, $\rho$ denotes the fluid density and $\mu$ denotes the fluid viscosity; we set both to $1$ for simplicity. The final time is denoted by $T_f$.  Since no time derivative of the pressure appears in the system, it is a DAE. The \textit{index} of a DAE is defined as the number of analytical differentiations needed (along with algebraic manipulations) to convert the DAE into an explicit system of ODEs~\cite[Section VII.1]{wanner1996solving}.  Here, since the constraint equation is of the form $-\nabla\cdot\textbf{u} = 0$, this is an index-two DAE, since one differentiation of the constraint (and applying the divergence to~\eqref{subeqn:velocity}) allows us to explicitly solve for $p$ in terms of $\textbf{u}$, and a second gives an ODE for $p$.  For index-two DAEs, the order of accuracy of a Runge–Kutta time-discretization is limited to the stage order of the scheme.

For the spatial discretization of \eqref{time stokes}, we use the mixed finite-element framework, considering the stable Taylor--Hood discretization on simplices \cite{elman}. Let $\mathcal{V} = \textbf{H}_{0}^{1}(\Omega)$, where $\textbf{H}_{0}^{1}(\Omega) =\{\textbf{v}\in \textbf{H}^{1}(\Omega): \textbf{u}=0 \text{ on }  \partial \Omega \},$ and $\mathcal{W}=L^2_0(\Omega)$ (the space of zero-mean functions in $L^2(\Omega)$), and consider a weak solution of \eqref{time stokes} that is (at least) once continuously differentiable in time and such that for every $t\in (0,T_f)$, $\textbf{u}(\cdot, t)\in \mathcal{V}$ and $p(\cdot, t)\in\mathcal{W}$.  Multiplying the time-dependent equation by $\textbf{v}\in\mathcal{V}$ and the divergence constraint by $q\in\mathcal{W}$ and integrating by parts, we get the weak form
\begin{alignat*}{3}
  \langle \vect{u}_t,\vect{v} \rangle + \langle \nabla \vect{u},\nabla\vect{v}\rangle - \langle p,\nabla\cdot\vect{v}\rangle & = \langle \vect{f},\vect{v}\rangle, &&\quad \forall\vect{v}\in\mathcal{V}, \\
  -\langle q,\nabla\cdot\vect{u}\rangle & = 0, &&\quad \forall q\in\mathcal{W},
\end{alignat*}
where the inner-product notation, $\langle\cdot,\cdot\rangle$, denotes integration in space but not time.
The finite-element discretization is realized by constructing a triangulation, $\tau_h$, of $\Omega$, and approximating $\vect{u}$ and $p$ in piecewise polynomial spaces defined over $\tau_h$.  Here, we use standard continuous Lagrange finite-element spaces, defining
\[
P_k(\Omega,\tau_h) = \left\{ u\in C^0(\Omega) : \forall T\in \tau_h, \left.u\right|_T(\vect{x}) \text{ is a polynomial of degree no more than }k\right\}.
\]
We consider the standard stable Taylor--Hood discretization, with $\mathcal{V}_h = \left(P_2(\Omega,\tau_h)\right)^d\cap \mathcal{V}$ and $\mathcal{W}_h = P_1(\Omega,\tau_h) \cap \mathcal{W}$ \cite{elman, taylor1973numerical}.  This leads to the semi-discretized weak form of finding $(\vect{u}(\cdot,t),p(\cdot,t))\in \mathcal{V}_h\times\mathcal{W}_h$ such that
\begin{alignat*}{3}
  \langle \vect{u}_t,\vect{v} \rangle + \langle \nabla \vect{u},\nabla\vect{v}\rangle - \langle p,\nabla\cdot\vect{v}\rangle & = \langle \vect{f},\vect{v}\rangle, &&\quad \forall\vect{v}\in\mathcal{V}_h, \\
  -\langle q,\nabla\cdot\vect{u}\rangle & = 0, &&\quad \forall q\in\mathcal{W}_h.
\end{alignat*}
Now writing $\vu(t)$ and $\vp(t)$ for the (time-dependent) coefficients of $\vect{u}(\mathbf{x},t)$ and $p(\mathbf{x},t)$ in the finite-element basis, we can write this as a coupled linear system of DAEs, as
\[
\begin{bmatrix} M\vu_t \\ 0 \end{bmatrix}
+ \begin{bmatrix} K & B \\ B^T & 0 \end{bmatrix}\begin{bmatrix} \vu \\ \vp \end{bmatrix} = \begin{bmatrix} M\vf \\ 0 \end{bmatrix},
\]
where $\vf$ is the vector of coefficients of the interpolant of $\vect{f}$ in $\mathcal{V}_h$.  Here, $M$ and $K$ are the $\left(P_2(\Omega,\tau_h)\right)^d$ mass and stiffness matrices, respectively, while $B$ is the weak gradient operator mapping from $\mathcal{W}_h$ into $\mathcal{V}_h$.

It is this system of equations that we discretize using Runge--Kutta methods.  As the system is a set of DAEs, and not ODEs, we cannot directly apply \eqref{general RK}, but use its DAE analogue \cite{wanner1996solving}, writing
\begin{alignat*}{2}
  & \vu_i^n = \vu^{n} + \Delta t \sum_{j=1}^r a_{ij}\vk_j^{(\vect{u})},
  &&
  \vp_i^n = \vp^{n} + \Delta t \sum_{j=1}^r a_{ij}\vk_j^{(p)}, \\
  & M\vk_i^{(\vect{u})} + K\vu_i^n + B\vp_i^n = M\vf_i^n, \qquad
  &&
  B^T\vu_i^n = 0, \\
  & \vu^{n+1} = \vu^n + \Delta t \sum_{j=1}^r b_{j}\vk_j^{(\vect{u})},
  &&
   \vp^{n+1} = \vp^n + \Delta t \sum_{j=1}^r b_{j}\vk_j^{p},
\end{alignat*}
where $\vf_i^n$ is the interpolant of $\vect{f}$ in $\mathcal{V}_h$ at time $t^n+c_i\Delta t$, $\vu_i^n$ and $\vp_i^n$ are the approximations of $\vu$ and $\vp$ at time $t^n + c_i\Delta t$, and $\vk_i^{(\vect{u})}$ and $\vk_i^{(p)}$ are the RK stages for which we solve.  Rewriting the equations for $\vk_i^{(\vect{u})}$ and $\vk_i^{(p)}$, we have
\begin{align*}
  M\vk_i^{(\vect{u})} + K\left(\vu^{n} + \Delta t \sum_{j=1}^r a_{ij}\vk_j^{(\vect{u})}\right) + B\left(\vp^{n} + \Delta t \sum_{j=1}^r a_{ij}\vk_j^{(p)}\right) & = M\vf_i^n \\
    B^T\left(\vu^{n} + \Delta t \sum_{j=1}^r a_{ij}\vk_j^{(\vect{u})}\right) & = 0
\end{align*}
or
\begin{align*}
  M\vk_i^{(\vect{u})} + \Delta t \sum_{j=1}^r a_{ij}\left(K\vk_j^{(\vect{u})} + B\vk_j^{(p)}\right) & = M\vf_i^n - K\vu^{n} - B\vp^n \\
  \Delta t \sum_{j=1}^r a_{ij}B^T\vk_j^{(\vect{u})} & = -B^T\vu^{n}
\end{align*}
for $1\leq i \leq r$.  The matrix on the left can easily be written in tensor-product form, leading to a concise description of the scheme as
 \begin{equation}
\left( \textbf{I}_{r}\otimes \begin{bmatrix}
M & 0\\
0 &0\\
\end{bmatrix}\\ + \Delta t A \otimes \begin{bmatrix}
K & B\\
B^{T} &0\\
\end{bmatrix}\\ \right)\vec{\textbf{k}}=\vec{\textbf{F}},
\label{time-depend. stokes system}
 \end{equation}
 where $\vec{\textbf{k}}$ is the vector of stages, ordered consecutively by stage index $i$, keeping the ordering of $(\vk_i^{(\vect{u})},\vk_i^{(p)})$ pairs together, and $\vec{\textbf{F}}$ is the corresponding vector of right-hand sides (including terms from timestep $n$).

 \subsection{Navier--Stokes Equations}

We next include the full inertial term, leading to the nonlinear incompressible Navier--Stokes equations,
\begin{subequations}
\label{navier-stokes}
\begin{alignat}{3}
\rho\left(\textbf{u}_{t}+ \textbf{u}\cdot \nabla \textbf{u}\right)-\mu \Delta \textbf{u} + \nabla p&=\vect{f}    && \text{  in } \Omega \times (0,T_{f}) , \\
-\nabla \cdot \textbf{u}&=0     && \text{  in } \Omega \times (0,T_{f}) ,\\
\textbf{u}&=0      && \text{  on } \partial\Omega  \times (0,T_{f}),\\
\textbf{u}(\textbf{x},0)&=\textbf{g}(\textbf{x})     &&\text{  on } \Omega \times \{ t=0 \}.
\end{alignat}
\end{subequations}
  We again take the density to be 1, but will allow the viscosity $\mu$ to be chosen differently, to consider problems at different Reynolds numbers. The additional term passes directly to the weak form, which we again discretize using a Taylor--Hood mixed finite-element discretization. The semi-discretized weak variational form of \eqref{navier-stokes} is to find $(\textbf{u}(\cdot,t),p(\cdot,t))\in \mathcal{V}_{h}\times \mathcal{W}_{h}$ such that
  \begin{equation}
\begin{aligned}
\langle\textbf{u}_{t},\textbf{v}\rangle+\langle\textbf{u}\cdot \nabla \textbf{u},\textbf{v}\rangle+\mu \langle\nabla \textbf{u}, \nabla \textbf{v}\rangle-\langle p,\nabla \cdot \textbf{v}\rangle&=0,\\
-\langle\nabla \cdot \textbf{u},q\rangle&=0,
\end{aligned}
 \label{semidiscrete navier}
 \end{equation}
  for all test functions $(\textbf{v},q)\in \mathcal{V}_{h}\times \mathcal{W}_{h}$.  Again writing $\vu(t)$ and $\vp(t)$ for the coefficients of $\vect{u}$ and $p$ in the finite-element basis, this leads to a nonlinear coupled system of DAEs, as
  \[
  \begin{bmatrix} M\vu_t \\ 0 \end{bmatrix}
  + \begin{bmatrix} N(\vu) \\ 0 \end{bmatrix}
  + \begin{bmatrix} K & B \\ B^T & 0 \end{bmatrix}\begin{bmatrix} \vu \\ \vp \end{bmatrix} = \begin{bmatrix} M\vf \\ 0 \end{bmatrix},
  \]
  where $N(\vu)$ represents the discretization of $\langle\textbf{u}\cdot \nabla \textbf{u},\textbf{v}\rangle$.  As above, accounting for this term in the RK stage equations leads to the nonlinear coupled system
\begin{align*}
  M\vk_i^{(\vect{u})} + N\left(\vu^{n} + \Delta t \sum_{j=1}^r a_{ij}\vk_j^{(\vect{u})} \right) + \Delta t \sum_{j=1}^r a_{ij}\left(K\vk_j^{(\vect{u})} + B\vk_j^{(p)}\right) & = M\vf_i^n - K\vu^{n} - B\vp^n, \\
  \Delta t \sum_{j=1}^r a_{ij}B^T\vk_j^{(\vect{u})} & = -B^T\vu^{n},
\end{align*}
for $1\leq i \leq r$.  This system is solved using Newton's method.

Denoting the nonlinear system as $F\left(\vbk^{n}\right) = 0$, a standard Newton approximation would be to solve
\[
F\left(\vbk^{n,\ell+1}\right)\approx F\left(\vbk^{n,\ell} \right)+  J\left(\vbk^{n,\ell}\right) \delta\vbk^{n,\ell} = 0,
\]
where $J\left(\vbk^{n,\ell}\right)$ is the Jacobian of the system at the current approximation, $\vbk^{n,\ell}$ and $\delta\vbk^{n,\ell} \coloneqq \vbk^{n,\ell+1} - \vbk^{n,\ell}$ is the Newton search direction.  Since we are timestepping, we use the computed solution at the previous time-step, $\vbk^{n-1}$, for the initial guess for the stage values at step $n$, $\vbk^{n,0}$.  In this work, we use the Eisenstat--Walker stopping criterion for the Krylov iteration to solve for $\delta\vbk^{n,\ell}$ \cite{eisenstat1996choosing}, requiring that
\[
\left\|F\left(\vbk^{n,\ell} \right)+  J\left(\vbk^{n,\ell}\right) \delta\vbk^{n,\ell}\right\| \leq \eta_\ell \left\|F\left(\vbk^{n,\ell} \right)\right\|,
\]
for every step, $\ell$, where $\eta_{\ell}\in [0,1)$ is updated for each nonlinear iteration based on convergence of the method. 
 
 \subsection{Magnetohydrodynamics}
 Finally, we consider the equations of single-fluid viscoresistive magnetohydrodynamics (MHD). In general, MHD models the flow of conducting fluids in the presence of an electromagnetic field. These models are nonlinear and contain strong coupling between the fluid velocity and the electromagnetic variables. We follow the MHD formulation presented in \cite{scottmhd}, 
\begin{subequations}
\label{MHD}
\begin{align}
 \textbf{u}_{t} + (\textbf{u}\cdot \nabla)\textbf{u}-\nabla \cdot ( \frac{2}{\text{Re}}\epsilon(\textbf{u}))+\nabla p-(\nabla \times \textbf{B})\times \textbf{B}&=\textbf{f}_{\textbf{u}}  \text{  in } \Omega \times (0,T_{f}) ,\\
\textbf{B}_{t}+\frac{1}{\text{Re}_{m}}\nabla \times \nabla \times \textbf{B} - \nabla \times(\textbf{u}\times\textbf{B})-\nabla \gamma&=\textbf{f}_{\textbf{B}}  \text{  in } \Omega \times (0,T_{f}) ,\\
 -\nabla \cdot \textbf{u}&=0  \text{  in } \Omega \times (0,T_{f}) ,\label{u const}\\
 \nabla \cdot \textbf{B}&=0  \text{  in } \Omega \times (0,T_{f}) ,\label{B const}\\
\textbf{u}&=0     \text{  on } \partial\Omega \times (0,T_{f}),\\
\textbf{B}\times \textbf{n}&=0     \text{  on } \partial\Omega \times (0,T_{f}),\\
\textbf{u}(\textbf{x},0)=\textbf{g}_{\textbf{u}}&(\textbf{x})      \text{  on } \Omega \times \{ t=0 \},\\
\textbf{B}(\textbf{x},0)=\textbf{g}_{\textbf{B}}&(\textbf{x})      \text{  on } \Omega \times \{ t=0 \},
\end{align}
\end{subequations}

where the four unknowns are the velocity vector, $\textbf{u}$, the pressure, $p$, the magnetic field, $\textbf{B}$, and the Lagrange multiplier, $\gamma$.  The Lagrange multiplier is used to enforce the solenoidal condition~\eqref{B const}, while the pressure is used to enforce the incompressibility condition \eqref{u const}. The strain-rate tensor is $\epsilon(\textbf{u})=\frac{1}{2}(\nabla \textbf{u}+\nabla \textbf{u}^{T})$, and the two dimensionless constants, $\text{Re}$ and $\text{Re}_{m}$, are the hydrodynamic Reynolds number and magnetic Reynolds number, respectively.  We consider this equation in both two- and three-dimensional domains, $\Omega$; in 2D, the curl and cross-product are defined by the natural extensions from two-dimensional vector fields to three-dimensional fields.
 
Here, for $\Omega \subset\mathbb{R}^d$, we take
\[
(\textbf{u}(\cdot,t),\textbf{B}(\cdot,t),p(\cdot,t),\gamma(\cdot,t))\in\textbf{H}_{0}^{1}(\Omega) \times \textbf{H}_{0}(\text{curl},\Omega)\times L_{0}^{2}(\Omega)\times H_{0}^{1}(\Omega),
\]
where
 \begin{align*}
\textbf{H}_{0}^{1}(\Omega) &=\{\textbf{v}\in \textbf{H}^{1}(\Omega): \textbf{u}=0 \text{ on }  \partial \Omega \}, \\
\textbf{H}_{0}(\text{curl},\Omega )&=\{\textbf{c}\in \textbf{L}^{2}(\Omega): \nabla \times \textbf{c} \in \textbf{L}^{2}(\Omega), \textbf{n}\times\textbf{c}=0 \text{ on }\partial \Omega \} ,\\
L_{0}^{2}(\Omega) &=\{q \in L^{2}(\Omega):\int_{\Omega}q \dx =0 \},\\
H_{0}^{1}(\Omega) &=\{s \in H^{1}(\Omega): s=0 \text{ on } \partial \Omega \},
 \end{align*}
and $\textbf{n}$ is the outward unit normal vector on $\partial \Omega$ \cite{scottmhd, Schneebeli}.  We discretize the fluid variables again with the Taylor--Hood discretization $\mathcal{V}_h\times\mathcal{W}_h \subset \textbf{H}_{0}^{1}(\Omega)\times L_{0}^2(\Omega)$, and use lowest-order N\'{e}d\'{e}lec elements for $\mathcal{C}_{h} \subset \textbf{H}_0(\text{curl},\Omega)$ and $\mathcal{S}_h = P_1(\Omega, \tau_h) \cap H^1_0(\Omega)$ for the Lagrange multiplier. Well-posedness (under small-data assumptions) of both the continuous and discrete formulations is shown in  \cite{Schneebeli}.  

Multiplying \eqref{MHD} by the test functions $(\textbf{v},\textbf{c},q,s)\in\mathcal{V}_{h}\times \mathcal{C}_{h}\times \mathcal{W}_{h}\times \mathcal{S}_{h}$ and integrating by parts, we get the semi-discretized weak variational form of finding $(\textbf{u}(\cdot,t),\textbf{B}(\cdot,t),p(\cdot,t),\gamma(\cdot,t))\in\mathcal{V}_{h}\times \mathcal{C}_{h}\times \mathcal{W}_{h}\times \mathcal{S}_{h}$ such that
\begin{equation}
 \begin{aligned}
\int_{\Omega} \textbf{u}_{t}\cdot \textbf{v} + \left((\textbf{u}\cdot \nabla)\textbf{u}\right)\cdot \textbf{v}-  \frac{2}{\text{Re}} \left(\epsilon(\textbf{u}):\epsilon(\textbf{v}) \right)-p\nabla\cdot \textbf{v} -\left((\nabla \times \textbf{B})\times \textbf{B}\right)\cdot \textbf{v} \dx \\  =  \int_{\Omega} \textbf{f}_{\textbf{u}}\cdot \textbf{v} \dx,\\
\int_{\Omega} \textbf{B}_{t}\cdot \textbf{c}+\frac{1}{\text{Re}_{m}}(\nabla \times \textbf{B})\cdot( \nabla \times \textbf{c}) - (\textbf{u}\times \textbf{B}) \cdot ( \nabla\times \textbf{c})-(\nabla \gamma )\cdot \textbf{c} \dx = \int_{\Omega}\textbf{f}_{\textbf{B}}\cdot \textbf{c} \dx,\\
 -\int_{\Omega} (\nabla \cdot \textbf{u}) q \dx =0,\\
 -\int_{\Omega} \textbf{B}\cdot (\nabla s) \dx =0,
  \end{aligned}
\label{MHD semi-discrete}
\end{equation}
for all $(\textbf{v},\textbf{c},q,s)\in\mathcal{V}_{h}\times \mathcal{C}_{h}\times \mathcal{W}_{h}\times \mathcal{S}_{h}$.  The corresponding IRK discretization is derived from this semi-discretized form as described above, and solved in the same Newton--Krylov--multigrid manner, using the Eisenstat--Walker stopping criterion for the Krylov iteration.  We note that linear solvers for this spatial discretization using  BDF2 in time was the subject of \cite{scottmhd}; given the Dahlquist barrier, and the driving need for L-stability, multistage schemes, such as IRK, are needed to achieve higher-order time integration for this problem.

\section{Monolithic multigrid for fluid problems}
\label{sec:monolithicMG}

As described above, we use a Newton--Krylov--multigrid framework for solving the non-linear systems of equations resulting from spatial and temporal discretization of the models in~\Cref{sec:spatial_disc} (noting that the Newton linearization is trivial in the case of the linear Stokes equations).  Since the systems are nonsymmetric, we use FGMRES \cite{saad1993flexible} as the Krylov method, and seek to effectively precondition it.  For IRK discretizations of scalar PDEs, such as the heat equation, block-diagonal preconditioning of the stage-coupled linear systems is known to be effective~\cite{mardal2007order}.  While block-diagonal preconditioning has also been developed for fluid models discretized using BDF-like methods~\cite{elman,Wathen_etal_2017a}, we leave extension of this approach to IRK discretizations for future work.  Instead, we follow the approach of~\cite{vanka1986block,van2005multigrid}, and develop a monolithic multigrid preconditioner that makes use of an overlapping additive Schwarz relaxation that can be viewed as the extension of Vanka relaxation to the IRK case.  Compared to the use of block-structured preconditioners, this offers the advantage of not needing to explicitly approximate Schur complements in the stage-coupled IRK linearizations (which may depend on properties of the specific IRK scheme chosen, for example).  We note that FGMRES and classical right-preconditioned GMRES solve the same underlying optimization problem for the approximation in the same Krylov space, but that the underlying algorithms have important differences, with FGMRES requiring extra vector storage (to store both the Arnoldi vectors and their preconditioned counterparts) but right-preconditioned GMRES requiring an extra application of the preconditioner once the solution to the underlying Hessenberg system has been found.  Thus, we choose to use FMGRES, instead of classical right-preconditioned GMRES, primarily because the cost of application of our preconditioners is non-trivial, but we are not memory-bound on the parallel machine used in the numerical results.  Thus, the extra vector storage of FGMRES is an attractive trade-off over the extra preconditioner application required by standard right-preconditioned GMRES. An auxiliary advantage is that FGMRES allows the use of GMRES inside inner iterations (such as the relaxation).

\subsection{Coarse-grid correction and transfer operators}
In the numerical results that follow, we consider hierarchies of grids generated by taking uniform refinements of a given coarsest grid.  To map functions from a coarse mesh to its refinement, we use canonical finite-element interpolation operators for each field in the discretization.  For the single-stage case, for both the time-dependent Stokes and Navier--Stokes problems, interpolation takes the form
$$
P=\begin{bmatrix}
P_{\textbf{u}} & \\
& P_{p}\\
\end{bmatrix},
$$
where $P_{\textbf{u}}$ and $P_{p}$ represent the interpolation operators for the $P_{2}$ and $P_{1}$ finite-element spaces, respectively. In the MHD case, we introduce finite-element interpolation operators $P_{\textbf{B}}$ for the lowest-order  N\'{e}d\'{e}lec space and $P_{\gamma}$ for the $P_{1}$ space (with suitable boundary conditions for $\gamma$), following \cite{scottmhd}, making the interpolation operator for the single-stage case
 $$
P=\begin{bmatrix}
P_{\textbf{u}} & & &\\
& P_{\textbf{B}} & &\\
& &P_{p}&\\
&&&P_{\gamma}\\
\end{bmatrix}.
$$ 
For multistage IRK discretizations, the interpolation operator is defined as $I_r\otimes P$, creating a block-diagonal interpolation operator that applies the finite-element interpolation in $P$ to each stage independently.  We use the transpose of interpolation as the restriction operator.

We use rediscretization to define the coarse-grid operators, noting that this is equivalent to Galerkin coarsening in the finite-element case (if compatible quadrature rules are used to assemble on the fine and coarse grids). The coarsest-grid systems are solved directly, using MUMPS~\cite{MUMPS_refs}.

\subsection{Vanka Relaxation}

Vanka relaxation was first introduced for the time-steady MAC-scheme discretization of the Navier--Stokes equations~\cite{vanka1985block}, but it has been extensively used in many more general settings in recent years~\cite{farrell2020local, DGstokes, adler2016monolithic, scottmhd, pcpatch}.
Broadly defined, Vanka relaxation schemes are overlapping Schwarz (domain decomposition) methods used as relaxation for a multigrid algorithm. In order to achieve the expected cost of a multigrid relaxation scheme, the subdomain problems are generally quite small, on the order of 10s-100s of DoFs.  Historically, the most common approaches were multiplicative in nature; however, we follow the recent trend towards additive schemes \cite{farrell2020local, scottmhd, pcpatch} that are naturally parallelizable.

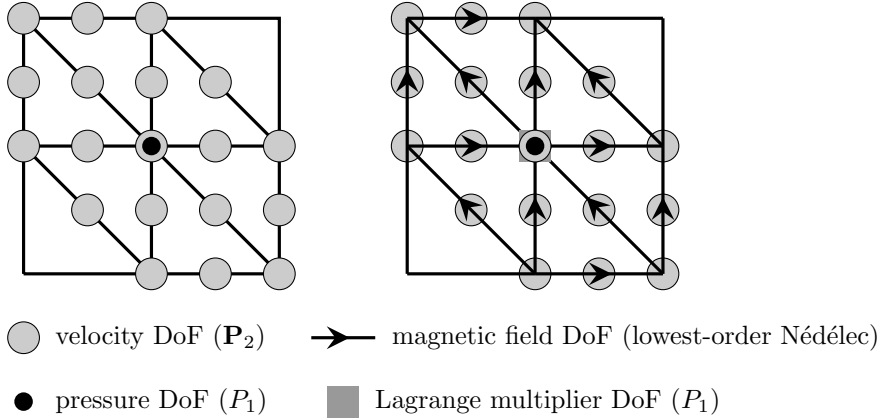
\begin{figure}
        \centering
        \begin{tikzpicture}[scale=0.85]
        %\draw[help lines] (0,0) grid (13,4);

        %% Fluid Patch
        \draw[very thick] (0,0) -- (4,0) -- (4, 4) -- (0, 4) -- (0,0);
        \draw[very thick] (0,2) -- (4,2);
        \draw[very thick] (2,0) -- (2,4);
        \draw[very thick] (4,0) -- (0,4);
        \draw[very thick] (2,0) -- (0,2);
        \draw[very thick] (4,2) -- (2,4);
        \foreach \x in {1,...,3}{
                \foreach \y in {1,...,3}{
                        \node[circle,draw,fill=gray!40,inner sep=0pt,minimum size=12pt] at (\x,\y) {};
                }
        }
        \node[circle,draw,fill=gray!40,inner sep=0pt,minimum size=12pt] at (2,0) {};
        \node[circle,draw,fill=gray!40,inner sep=0pt,minimum size=12pt] at (2,4) {};
        \node[circle,draw,fill=gray!40,inner sep=0pt,minimum size=12pt] at (0,2) {};
        \node[circle,draw,fill=gray!40,inner sep=0pt,minimum size=12pt] at (4,2) {};
        \node[circle,draw,fill=gray!40,inner sep=0pt,minimum size=12pt] at (4,0) {};
        \node[circle,draw,fill=gray!40,inner sep=0pt,minimum size=12pt] at (3,0) {};
        \node[circle,draw,fill=gray!40,inner sep=0pt,minimum size=12pt] at (4,1) {};
        \node[circle,draw,fill=gray!40,inner sep=0pt,minimum size=12pt] at (0,3) {};
        \node[circle,draw,fill=gray!40,inner sep=0pt,minimum size=12pt] at (0,4) {};
        \node[circle,draw,fill=gray!40,inner sep=0pt,minimum size=12pt] at (1,4) {};
        \node[circle,fill=black,inner sep=0pt,minimum size=7pt] at (2,2) {};
        \node at (1.6,1.6) {};

        \tikzset{->-/.style={decoration={markings,mark=at position .6 with {\arrow[scale=1.7]{stealth}}},postaction={decorate}}}

        %% Combined Patch
        \foreach \x in {7,...,9}{
                \foreach \y in {1,...,3}{
                        \node[circle,draw,fill=gray!40,inner sep=0pt,minimum size=12pt] at (\x,\y) {};
                }
        }
        \node[circle,draw,fill=gray!40,inner sep=0pt,minimum size=12pt] at (8,0) {};
        \node[circle,draw,fill=gray!40,inner sep=0pt,minimum size=12pt] at (8,4) {};
        \node[circle,draw,fill=gray!40,inner sep=0pt,minimum size=12pt] at (6,2) {};
        \node[circle,draw,fill=gray!40,inner sep=0pt,minimum size=12pt] at (10,2) {};
        \node[circle,draw,fill=gray!40,inner sep=0pt,minimum size=12pt] at (10,0) {};
        \node[circle,draw,fill=gray!40,inner sep=0pt,minimum size=12pt] at (9,0) {};
        \node[circle,draw,fill=gray!40,inner sep=0pt,minimum size=12pt] at (10,1) {};
        \node[circle,draw,fill=gray!40,inner sep=0pt,minimum size=12pt] at (6,3) {};
        \node[circle,draw,fill=gray!40,inner sep=0pt,minimum size=12pt] at (6,4) {};
        \node[circle,draw,fill=gray!40,inner sep=0pt,minimum size=12pt] at (7,4) {};
        \draw[very thick] (6,0) -- (8,0);
        \draw[very thick,->- ] (8,0) -- (10,0);
        \draw[very thick,->-] (6,2) -- (8,2);
        \draw[very thick,->-] (8,2) -- (10,2);
        \draw[very thick,->-] (6,4) -- (8,4);
        \draw[very thick] (8,4) -- (10,4);
        \draw[very thick] (6,0) -- (6,2);
        \draw[very thick,->-] (6,2) -- (6,4);
        \draw[very thick,->-] (8,0) -- (8,2);
        \draw[very thick,->-] (8,2) -- (8,4);
        \draw[very thick,->-] (10,0) -- (10,2);
        \draw[very thick] (10,2) -- (10,4);
        \draw[very thick,->-] (8,0) -- (6,2);
        \draw[very thick,->-] (10,0) -- (8,2);
        \draw[very thick,->-] (8,2) -- (6,4);
        \draw[very thick,->-] (10,2) -- (8,4);
        \node[rectangle,fill=white,gray!80,inner sep=0pt,minimum size=12pt] at (8,2) {};
        \node at (8.4,2.4) {};
        \node[circle,draw,fill=gray!40,inner sep=0pt,minimum size=12pt] at (8,2) {};
        \node[circle,fill=black,inner sep=0pt,minimum size=7pt] at (8,2) {};
        \node at (7.6,1.6) {} ;

% Legend
        \node[circle,draw,fill=gray!40,inner sep=0pt,minimum size=12pt] at (0,-1) {};
        \node[right=3mm] at (0,-1) {velocity DoF ($\mathbf{P}_2$)};
        \draw[very thick,->-] (4.5,-1) -- (5.5,-1) {};
        \node[right=1mm] at (5.5,-1) {magnetic field DoF (lowest-order N\'{e}d\'{e}lec)};
        \node[circle,fill=black,inner sep=0pt,minimum size=7pt] at (0,-2) {};
        \node[right=3mm] at (0,-2) {pressure DoF ($P_1$)};
        \node[rectangle,fill=white,gray!80,inner sep=0pt,minimum size=12pt] at (5,-2) {};
        \node[right=3mm] at (5,-2) {Lagrange multiplier DoF ($P_1$)};

        \end{tikzpicture}
        \caption{Left: Vanka patch for the Stokes and Navier--Stokes equations, consisting of $P_{2}$ velocity DoFs, and one $P_{1}$ pressure DoF. Right: Vanka patch for the MHD equations, consisting of $P_{2}$ velocity DoFs, lowest-order N\'{e}d\'{e}lec DoFs for the magnetic field, one $P_{1}$ Lagrange multiplier DoF and one $P_{1}$ pressure DoF. }\label{fig:vanka_patches}
\end{figure}

To specify the details of relaxation, we now describe how the Schwarz subdomains (commonly referred to as the Vanka ``blocks'' or ``patches'') are constructed from the underlying mesh on any given level of the multigrid hierarchy.  In this work, we follow the topological construction described in \cite{pcpatch}.  In particular, we form a Vanka patch for each vertex in the mesh, which consists of all degrees of freedom associated with the closure of the cells adjacent to the vertex.  As is typical in Vanka relaxation, we exclude all degrees of freedom associated with $P_1$ constraints (the pressure and Lagrange multiplier in our systems) from the patch, except for those located at the vertex around which the patch is formed.  For the models considered here, this results in patches like those shown in~\Cref{fig:vanka_patches} for regular two-dimensional grids, with a single pressure degree of freedom and all velocity DoFs on all elements adjacent to the node.  When used in an IRK discretization, these patches include all stage degrees of freedom.  For MHD, we note that the patch shown at right of~\Cref{fig:vanka_patches} coincides topologically with the \textit{coupled Vanka} approach for the BDF2 discretization considered in~\cite{scottmhd}, but the patches used here contain more degrees of freedom than those used in~\cite{scottmhd}, due to inclusion of all stages in the IRK discretization.

Denoting the set of DoFs in the $i^{th}$ Vanka patch by $\mathcal{S}_{i}$, we have (by construction) that every degree of freedom in the domain is contained in at least one patch: $\mathcal{S}=\bigcup^{N}_{i=1} \mathcal{S}_{i}$, where $N$ is  the total number of patches and $\mathcal{S}$ is the complete set of DoFs for the problem. Denoting $R_{i}$ as a ``restriction'' operator that maps global DoFs to those in patch $\mathcal{S}_{i}$, we can write a single iteration of a weighted stationary iteration as
 $$\vec{\mathbf{k}}\leftarrow \vec{\mathbf{k}}+\omega\sum _{i=1} ^{N} R_{i}^{T}(R_{i}JR_{i}^{T})^{-1}R_{i}(\vec{\mathbf{F}}-J\vec{\mathbf{k}}),$$
where $J\vec{\mathbf{k}} = \vec{\mathbf{F}}$ is the linear system to be solved, and $R_{i}JR_{i}^{T}$ is the restriction of $J$ to the DoFs in patch $\mathcal{S}_{i}$. 
In practice, we use several steps of a Vanka-preconditioned Chebyshev or GMRES iteration as the relaxation scheme for our problems, with the endpoints of the interval defining the associated Chebyshev polynomials tuned by hand.

\subsection{Implementation}
The numerical results below are produced using Firedrake~\cite{rathgeber2016firedrake} for the spatial finite-element discretization and Irksome~\cite{Irksome} for the temporal discretization.  Linear and nonlinear solvers are implemented in PETSc~\cite{balay2018petsc}, taking advantage of the close integration between discretizations and solvers provided by this combination~\cite{kirby2018solver}.  The Vanka relaxation is implemented through PCPATCH~\cite{pcpatch}.  For reproducibility, the codes used to generate the numerical results and the major components of Firedrake, Irksome, and PETSc needed, have been archived on Zenodo~\cite{zenodo/Firedrake-20220210.0}.
We emphasize that all aspects of the discretization and solver software are chosen to be naturally parallelizable.  The coarsest mesh in each hierarchy is distributed across the available parallel cores, and then refined in parallel.  For the two-dimensional problems below, after each refinement, the mesh is redistributed to better balance parallel work, but this was not done for the 3D example, due to software limitations.  While not rebalancing the meshes leads to a small load imbalance on the finer grids in the hierarchies, this was not seen to lead to significant loss of performance in the weak scaling tests reported below.  To account for the need to compute residuals for each DoF in each Vanka patch, a node-distance-2 halo is included in the parallel mesh distribution, to allow the relaxation scheme to be performed in parallel without additional communication~\cite{dmplex}.

\section{Numerical Results}
\label{sec:numerics}

For the numerical results in this paper, we focus on 3 families of IRK methods:  Gauss (also known as Gauss--Legendre), LobattoIIIC, and RadauIIA.  We note that LobattoIIIC and RadauIIA are both L-stable and A-stable, while Gauss is A-stable, but not L-stable.  All are fully implicit schemes, with stage order equal to the number of stages.  For an $r$-stage method, Gauss schemes have order $2r$, while RadauIIA have order $2r-1$ and LobattoIIIC have order $2r-2$.  We consider both 2- and 3-stage schemes here, with standard Butcher tableaux~\cite{wiki_IRKlist, wanner1996solving}.

 We present results for four separate test cases: a simple two-dimensional time-dependent Stokes model, two-dimensional Navier--Stokes flow past a cylinder, and two MHD examples, a two-dimensional island-coalescence problem and a three-dimensional lid-driven cavity model.
 All results presented in this paper were computed on the Compute Canada cluster, \textit{Niagara}, consisting of 2,024 nodes, each with 40 2.4 GHz Intel Skylake cores and 202GB of RAM, connected using a 100Gb/s EDR Dragonfly+ network.
 
\subsection{Two-dimensional time-dependent Stokes}

We consider a method of manufactured solutions test case, solving~\eqref{time stokes} on the two-dimensional unit square, $\Omega =(0,1)^{2}$. The forcing function $\textbf{f}$ and boundary conditions are chosen so that the exact solution is $$\textbf{u}=\begin{bmatrix}\sin(\pi x)\cos(\pi y)e^{-2t\pi^{2}}\\ -\cos(\pi x)\sin(\pi y)e^{-2t\pi^{2}} \end{bmatrix}, \text{ and }p=0.$$  For this example, we construct a coarsest grid by creating a uniform $8\times 8$ quadrilateral mesh of the unit square, then cut each quadrilateral cell into 4 triangles, adding a vertex at the center of the quadrilateral.  This mesh is then uniformly refined $\ell$ times; below, we present results for $\ell=5,6,7$, where the Taylor--Hood discretization of the Stokes equations results in about 1.1 million DoFs per stage for $\ell = 5$ up to about 19 million DoFs per stage for $\ell=7$.  The initial condition is chosen by interpolating the exact solution into the finite-element space at $t=0$, and we integrate up to time $T_f = 0.5$, with timestep $\Delta t= T_f/N$ for $N=2^{\ell +3}$. To our knowledge, there are no rigorous stopping tolerances that guarantee discretization-error level accuracy for these systems; we use a hand-tuned stopping tolerance, where we require the absolute value of the $\ell_2$ norm of the residual of the system to be reduced below $10^{-2}\times N^{-3}$ at each timestep, or a relative reduction in this norm by $10^{-8}$.  Based on preliminary experiments, we accelerate the relaxation process using Chebyshev polynomials of the first kind on the interval $[2,8]$ and employ 2 pre- and post-relaxation sweeps. Proper choice of the Chebyshev interval is critical to achieving scalable performance.  For two-dimensional problems with geometric coarsening by a factor of two (as used here), a reasonable strategy is to estimate the largest eigenvalue, $\lambda$, of the relaxation-preconditioned matrix (e.g., using Ritz values from preconditioned GMRES) and choose the interval to be $[\lambda/4,\lambda]$.  Here, we started from similar estimates, but hand-tuned the intervals to optimize performance.

\Cref{table: 2d stokes results} presents a weak scaling study for this problem, for both two- and three-stage methods. For the two-stage methods, we use 10 cores on 1 node for $\ell=5$, 40 cores on 1 node for $\ell=6$ and 160 cores on 4 nodes for $\ell=7$. For the three-stage methods, we increase core counts by 50\%, to account for the increased number of degrees of freedom in the resulting linear systems, using 15 cores on 1 node for $\ell=5$, 60 cores on 2 nodes for $\ell=6$ and 240 cores on 6 nodes for $\ell=7$. We report the relative $L_2$ error in the velocity and the absolute $L_2$ error in the pressure approximation at the final time, as well as the average number of linear iterations to achieve convergence over all timesteps and the total computational time needed in minutes. In the final column of \Cref{table: 2d stokes results}, we report the average wall-clock time per Krylov iteration (t/K) in seconds.

\begin{table}
\centering
\begin{tabular}{c|c|c|c|c|c|c} 
\toprule
 & $\ell$ & velocity error & pressure error&iterations&time& t/K\\
\midrule
\multirow{3}{*}{Gauss(2)} &5& $1.786 \times 10^{-2}$  &	$2.328 \times 10^{-2}$ &	9.85 &	57.47 & 1.401\\
&6& $3.155 \times 10^{-3}$ &	$1.162 \times 10^{-2}$ &	13.39 &	202.98&1.779\\
&7& $5.575 \times 10^{-4}$ &	$6.271 \times 10^{-3}$ &	19.14 &	725.02&2.285\\   
\midrule
\multirow{3}{*}{RadauIIA(2)}  & 5 & $3.380 \times 10^{-6}$ & $6.327 \times 10^{-9}$ & 8.70 & 56.13&1.570 \\ 
 				& 6 & $4.297 \times 10^{-7}$ &	$2.795 \times 10^{-9}$ &	10.99 &	171.08&1.837\\
				& 7 & $4.971 \times 10^{-8}$ &	$1.028 \times 10^{-9}$ &	14.22 &	575.77 &2.391\\                                   
\midrule
\multirow{3}{*}{LobattoIIIC(2)} & 5 & $9.823 \times 10^{-4}$ & $4.594 \times 10^{-7}$ & 9.23  & 54.13 &1.912\\ 
                                    & 6 & $2.495 \times 10^{-4}$ & $1.479 \times 10^{-7}$ & 11.71 & 181.88 &1.833\\ 
                                    & 7 & $6.289 \times 10^{-5}$ & $4.823 \times 10^{-8}$ & 15.25 & 618.77 &2.393 \\  
\midrule
\midrule
\multirow{3}{*}{Gauss(3)}  &5&$9.098 \times 10^{-5}$ & $1.481 \times 10^{-4}$&	13.38&	109.39&1.929\\
&6&$1.839 \times 10^{-5}$ & $6.802 \times 10^{-4}$&	24.51&	483.75&2.323\\
                         & 7 & $3.293 \times 10^{-6}$ & $3.235 \times 10^{-4}$  & 25.97  &  1332.40&3.249 \\  
\midrule
\multirow{3}{*}{RadauIIA(3)} &5 &$6.151 \times 10^{-7}$ &	$1.298 \times 10^{-9}$ &	9.32 &	82.29&2.083\\
 &6 &$1.122 \times 10^{-7}$ &	$1.310 \times 10^{-9}$ &	12.40 &	267.78&2.566\\
 &7 &$1.300 \times 10^{-8}$ &	$1.718 \times 10^{-10}$ &	13.91 &	877.99&3.786\\  
\midrule
\multirow{3}{*}{LobattoIIIC(3)}  &5& $6.019 \times 10^{-7}$ &	$2.728 \times 10^{-9}$ &	9.61 &	85.25 &2.120\\
&6& $1.083 \times 10^{-7}$ &	$1.966\times 10^{-9}$ &	12.91 &	248.95&2.293\\
&7& $1.271 \times 10^{-8}$ &	$3.797 \times 10^{-10}$ &	14.91 &	938.06 &3.770\\ 
\bottomrule
\end{tabular}
\caption{Numerical results for two-dimensional Stokes model problem with two- and three-stage IRK schemes.  Relative $L_2$ errors in velocity and absolute $L_2$ errors for pressure are reported, along with average number of linear solver iterations per time-step, total wall-clock time-to-solution in minutes, and time per Krylov iteration in seconds, for refinement levels $\ell=5,6,7$.}
\label{table: 2d stokes results}
\end{table}

%[Personal note] To find the last column in above: Convert total run time on Niagara to seconds  then divide by (number of timesteps*its column) 

Table~\ref{table: orders of convergence 2d stokes results} summarizes rates of convergence for the results shown in \Cref{table: 2d stokes results}.
We observe at least second-order convergence in the velocity error for all three IRK schemes; however, we notice much larger errors for the Gauss results than for the other two schemes.  We note that the stopping tolerance decreases by a factor of 8 with each refinement, so that the slight increase in averaged iterations to convergence is not overly surprising, and seems to remain bounded at reasonable levels.  Nonetheless, the factor four increase in the number of cores with each refinement is insufficient to lead to ideal time scaling (which would be to double with each refinement, due to the doubling of the number of time-steps).

\begin{table}
\centering
\begin{tabular}{c|c|c|c|c|c|c|c} 
\toprule
            &  & \multicolumn{2}{|c|}{Gauss}& \multicolumn{2}{|c|}{RadauIIA}& \multicolumn{2}{|c}{LobattoIIIC}\\
\midrule
&&$u$&$p$ &$u$&$p$ &$u$&$p$\\
\midrule

\multirow{2}{*}{two-stage} & $\log_2{e_5}/{e_6}$ & 2.5& 1.0 &3.0 & 1.2 & 2.0 & 1.6  \\ 
                          \cline{2-8}				
 				& $\log_2{e_6}/{e_7}$ & 2.5 & 0.9 & 3.1 & 1.4& 2.0 & 1.6\\
\midrule				
\midrule
\multirow{2}{*}{three-stage}& $\log_2{e_5}/{e_6}$ &2.4& -2.3& 2.5& 0.0&2.5 & 0.5   \\ 
                          \cline{2-8}				
 				& $\log_2{e_6}/{e_7}$ & 2.5 & 1.1 & 3.1 & 2.9 & 3.1 & 2.4\\ 

\bottomrule
\end{tabular}
\caption{Rates of convergence in velocity and pressure for data in Table~\ref{table: 2d stokes results} with two- and three-stage IRK schemes for refinement levels $\ell=5,6,7$.  Here, $e_\ell$ denotes the error in a quantity on refinement level $\ell$.}
\label{table: orders of convergence 2d stokes results}
\end{table}

There are several contributing factors to the less-than-perfect scaling, beyond the simple increase in total number of Krylov iterations with refinement.  When going from $\ell=5$ to $\ell=6$ with the two-stage schemes, we increase the number of cores used for the calculation, but those cores remain on one physical node, leading to a saturation of the memory bandwidth available.  The same limitation occurs when going from $\ell=6$ to $\ell=7$ with the three-stage schemes, where we go from using 30 cores on each of 2 nodes to all 40 cores on 6 nodes.  While this could be avoided by using the same number of cores on more nodes of the parallel machine, such usage is impractical when a single node has sufficient memory for the $\ell=6$ problem with 2 IRK stages.  Furthermore, when going from $\ell=6$ to $\ell=7$, the (direct) coarsest-grid solve goes from being dominated by its computation to being dominated by its communication. Here, we clearly see another increase in the cost per linear iteration, especially in the three-stage methods where the time required increases by about 50\% for $\ell=7$. Improved performance would almost certainly be seen by duplicating the coarse-grid solve on each node, as considered in~\cite{9447889,10.1145/2929908.2929913,reisner2018scaling}.  We leave these performance enhancements for future work.

The experiment in \Cref{table: 2d stokes results} highlights convergence as we change both the spatial and temporal discretizations.  Here, however, we note that the temporal discretizations are higher order than the spatial, particularly for the 3-stage discretizations.  Thus, for comparison with this data, \Cref{table: 2d stokes results fixed dt} presents results with the same setup as \Cref{table: 2d stokes results}, but using a fixed timestep of $\Delta t=0.5/2^{8}$, to match results with $\ell = 5$ (noting that these results were run independently, so small differences in timings for $\ell = 5$ naturally arise).  We see that while quite reasonable convergence is observed in \Cref{table: 2d stokes results}, we observe significant stagnation in convergence here. Thus, even though the temporal discretizations are higher order, we still see substantial benefits to varying the timestep simultaneously with refinement of the spatial mesh. 

\begin{table}
\centering
\begin{tabular}{c|c|c|c|c|c|c} 
\toprule
 & $\ell$ & velocity error & pressure error&iterations&time& t/K\\
\midrule
\multirow{3}{*}{RadauIIA(2)}  & 5 & $3.380 \times 10^{-6}$ & $6.327 \times 10^{-9}$ & 8.70 & 56.37&1.560 \\ 
 				& 6 & $3.119 \times 10^{-6}$ &	$1.896 \times 10^{-9}$ &	12.18 &	96.77&0.940\\
				& 7 & $3.101 \times 10^{-6}$ &	$2.793 \times 10^{-9}$ &	18.04 &	165.42 &0.540\\                                   
\midrule
\multirow{3}{*}{RadauIIA(3)} &5 &$6.151 \times 10^{-7}$ &	$1.298 \times 10^{-9}$ &	9.32 &	82.07&2.087\\
 &6 &$3.677 \times 10^{-8}$ &	$3.288 \times 10^{-9}$ &	14.22 &	116.27&0.964\\
 &7 &$3.566 \times 10^{-8}$ &	$8.063 \times 10^{-10}$ &	17.30 &	238.18&0.809\\  
\midrule
\multirow{3}{*}{LobattoIIIC(3)}  &5& $6.019 \times 10^{-7}$ &	$2.728 \times 10^{-9}$ &	9.61 &	75.87 &1.850\\
&6& $6.066 \times 10^{-8}$ &	$3.076\times 10^{-9}$ &	14.60 &	107.87&0.866\\
&7& $5.611 \times 10^{-8}$ &	$1.136 \times 10^{-9}$ &	17.47 &	239.03 &0.805\\ 
\bottomrule
\end{tabular}
\caption{Results analogous to \Cref{table: 2d stokes results}, but with $\Delta t=0.5/2^{8}$.  Relative $L_2$ errors in velocity and absolute $L_2$ errors for pressure are reported, along with average number of linear solver iterations per time-step, total wall-clock time-to-solution in minutes, and time per Krylov iteration in seconds, for refinement levels $\ell=5,6,7$.}
\label{table: 2d stokes results fixed dt}
\end{table}

Finally, \Cref{table: 2d stokes DIRK results} presents comparison results for diagonal IRK schemes.  Here, we consider the two-stage second-order Pareschi-Russo (with parameter $1-\sqrt{2}/2$) \cite{PareschiRusso} and three-stage third-order Alexander \cite{Alexander} integrators, which are both L-stable.  While these results show some outperformance of the theoretical guarantees given by their stage order of one, they are also quite poor in comparison to the RadauIIA integrators of the same number of stages. In particular, comparing with the results in~\Cref{table: 2d stokes results}, we see that the errors achieved using DIRK(3) with $\ell=6$ are comparable to those achieved when using RadauIIA(2) with $\ell=5$, but that the latter calculation was achieved in about 60\% of the wall-clock time and on one-eighth of the number of cores (10 for RadauIIA(2) with $\ell=5$ vs. 80 for DIRK(3) with $\ell=6$).  Similarly, the errors for DIRK(3) with $\ell=7$ are slightly better than those achieved with RadauIIA(2) and $\ell=6$, and slightly worse than those achieved with RadauIIA(3) and $\ell=6$.  The two-stage Radau results, however, are achieved in just over 40\% of the wall-clock time, and on one-sixth the cores, while the three-stage Radau results are achieved in about two-thirds the wall-clock time, on one-fourth the cores.  These results highlight the added accuracy that can be gained using fully implicit RK methods over DIRK methods, and the added efficiency possible when using state-of-the-art linear solvers to achieve that accuracy.  For this reason, we focus on only the fully implicit RK schemes in the remainder of the paper.
\begin{table}
\centering
\begin{tabular}{c|c|c|c|c|c|c} 
\toprule
 & $\ell$ & velocity error & pressure error&iterations&time& t/K\\
\midrule
\multirow{3}{*}{DIRK(2)} &5 &$2.480 \times 10^{-5}$ &	$9.266 \times 10^{-8}$ &	6.42 &	35.20&1.312\\
 &6 &$6.198 \times 10^{-5}$ &	$3.254 \times 10^{-8}$ &	8.23 &	110.43&1.590\\
 &7 &$1.546 \times 10^{-5}$ &	$1.173 \times 10^{-8}$ &	10.07 &	213.76&1.835\\  
\midrule
\multirow{3}{*}{DIRK(3)}  & 5 & $1.106 \times 10^{-5}$ & $4.161 \times 10^{-8}$ & 6.99 & 35.28&1.240 \\ 
 				& 6 & $2.325 \times 10^{-6}$ &	$1.857 \times 10^{-9}$ &	9.51 &	$94.98^{*}$&1.12\\
				& 7 & $2.607 \times 10^{-7}$ &	$4.578 \times 10^{-9}$ &	11.17 &	415.26 &2.36\\                                   
\bottomrule
\end{tabular}

\caption{Numerical results for two-dimensional Stokes model problem with two- and three-stage DIRK schemes.  Relative $L_2$ errors in velocity and absolute $L_2$ errors for pressure are reported, along with average number of linear solver iterations per time-step, total wall-clock time-to-solution in minutes, and time per Krylov iteration in seconds, for refinement levels $\ell=5,6,7$.  Due to a change in the configuration of the machine, results for DIRK(3) at $\ell=6$ were run on 80 cores, instead of 60; all other results were run with same parallelism as in \Cref{table: 2d stokes results}.}
\label{table: 2d stokes DIRK results}
\end{table}

\subsection{Two-dimensional Navier--Stokes}

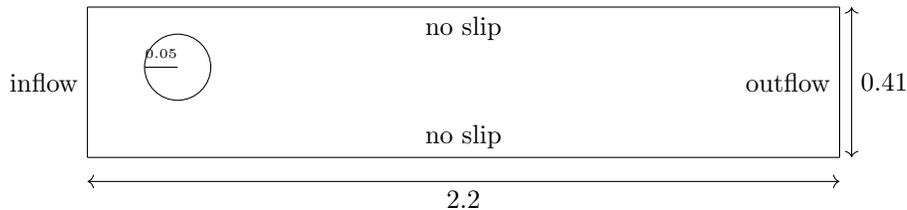
\begin{figure}
\begin{center}
  \begin{tikzpicture}[scale=0.8]

     \draw (-0.5,0) -- node[above] {no slip}(12,0) ;
     \draw (-0.5,0) -- node[left] {inflow}(-0.5,2.5) ;
     \draw (12,0) -- node[left] {outflow}(12,2.5) ;
     \draw (-0.5,2.5)-- node[below] {no slip}  (12,2.5) ;

     \draw (1,1.5) circle (0.55cm);
     \draw (0.45,1.5)--node[above]{\tiny 0.05}(1,1.5);
 
     \draw [<->](-0.5,-0.4)-- node[below] {2.2}  (12,-0.4) ;
     \draw [<->](12.2,0)-- node[right] {0.41}  (12.2,2.5) ;

   \end{tikzpicture}
\end{center}
\caption{Domain for Navier--Stokes flow past a cylinder.}
\label{navier_domain}
\end{figure}

We next consider two-dimensional Navier--Stokes flow past a cylinder, following the example given in~ \cite{john2004reference, Irksome, schafer1996benchmark}. Here, we consider the spatial domain $\Omega=(0,2.2)\times(0,0.41)\setminus B_{r}(0.2,0.2)$, where $B_r(0.2,0.2)$ is the disc of radius $r=0.05$ centred at $(0.2,0.2)$, shown in \Cref{navier_domain}. No-slip (zero-velocity) boundary conditions are imposed on the top and bottom boundaries of the rectangle and along the surface of the cylinder. Time-dependent inflow conditions are given on the left edge, prescribing
$$\textbf{u}(0,y,t)=\begin{bmatrix}\frac{4U(t)y(0.41-y)}{0.41^{2}}\\0\end{bmatrix},$$
 where $U(t)=1.5\sin\left(\frac{\pi t}{8}\right)$ is the mean inflow velocity. No-stress outflow is prescribed on the right boundary. The viscosity is set as $\mu =10^{-3}$, resulting in a Reynolds number of 100. The time step for these experiments is fixed as $\Delta t= \frac{1}{400}$, and we consider the final time $T_f = 8$.  As above, we discretize using Taylor--Hood elements in space and IRK in time.

For this problem, an unstructured coarsest grid with 972 triangular elements is used, chosen to refine the representation around the included cylinder.  Below, we report results for $3 \leq \ell \leq 6$, with discrete problem sizes for the Taylor--Hood discretization ranging from about 245 thousand DoFs per stage for $\ell=3$ to about 15.5 million DoFs per stage for $\ell=6$.  Details of the parallelization are provided in~\Cref{table: cores for Navier stokes}, where we again note that we have increased the number of cores for the 3-stage IRK methods by about 60\% over those for the 2-stage methods.  For this problem, we use a nonlinear stopping tolerance requiring the absolute $\ell_2$ norm of the nonlinear residual be below ${1}/{N^{3}}$ with $N= 2^{\ell +3}$, and use an Eisenstat--Walker inexact Newton scheme to determine the linear stopping tolerances for each nonlinear iteration.  Here, again 2 pre- and post-relaxation sweeps are used, with Chebyshev polynomials for relaxation taken over the interval $[1.5,8]$.
  
\begin{table}
\centering
\begin{tabular}{c|c|c|c|c} 
\toprule
  stages & $\ell$ & total DoFs & nodes & cores \\
\midrule
\multirow{4}{*}{2} & 3 & 489,656 & 1 & 6 \\
                   & 4 & 1,948,144 & 1 & 25 \\ 
                   & 5 & 7,771,616 & 4 & 100 \\ 
                   & 6 & 31,044,544 & 10 & 400 \\  
\midrule
\multirow{4}{*}{3} & 3 & 734,484 & 1 & 10 \\
                   & 4 & 2,922,216 & 1 & 40 \\ 
                   & 5 & 11,657,424 & 4 & 160 \\ 
                   & 6 & 46,566,816 & 16 & 640 \\                                    

\bottomrule
\end{tabular}
\caption{Total number of DoFs and number of nodes and cores used for the Navier--Stokes test problem with two- and three-stage IRK discretizations. }
\label{table: cores for Navier stokes}
\end{table} 

\begin{figure}
\centering
\includegraphics[width=0.45\textwidth]{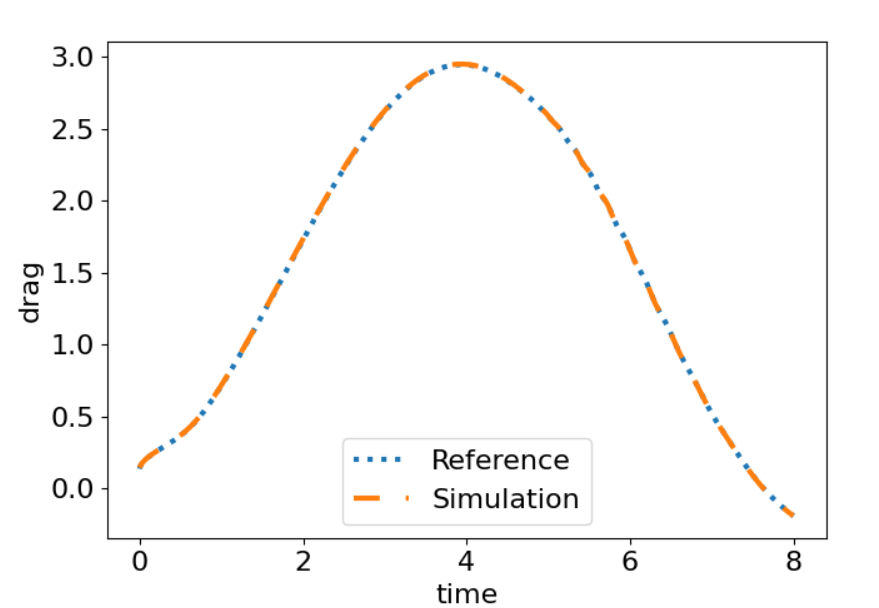}
%\subcaption{Drag}
\label{fig:subfig1}
\qquad
\includegraphics[width=0.45\textwidth]{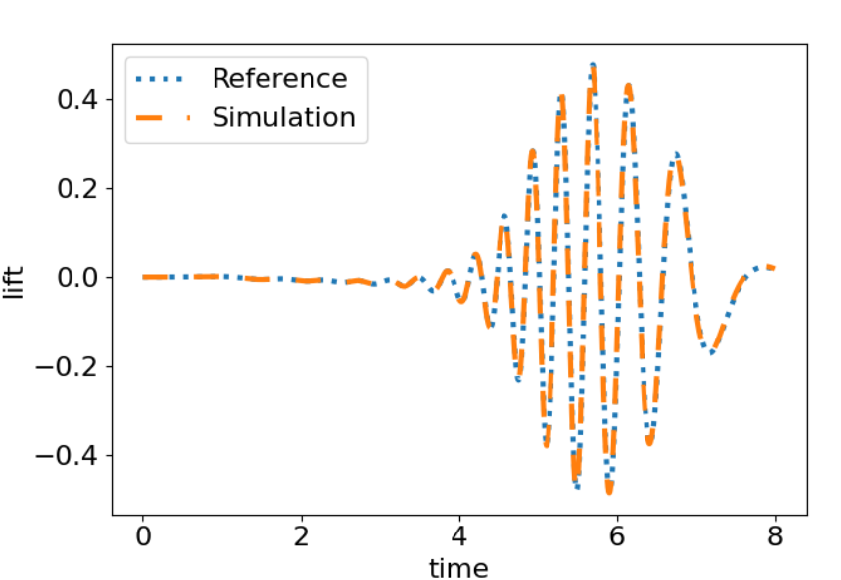}
%\subcaption{Lift}
\label{fig:subfig2}
\caption{Comparison of reference and drag (left) and lift (right) computed using LobattoIIIC(2) and $\ell=6$.}
\label{fig: NS figs}
\end{figure}

As no analytical solution is available in this case, we instead record the maximum drag and lift values computed over the simulations, for comparison with reference data~\cite{john2004reference, Irksome}.  \Cref{fig: NS figs} presents time-histories of these quantities for one simulation, showing excellent agreement with reference data. Results for other simulations with both RadauIIA and LobattoIIIC are visually similar.  \Cref{2d NS results} presents these values for both of these integrators, along with the average wall-clock time in minutes per time-step, and average nonlinear and linear solver iterations per time-step.  As before, we have decreasing solver tolerances as $\ell$ increases, so the small increases in iterations counts with refinement are expected.
\begin{table}
\centering
\begin{tabular}{c|c|c|c|c|c|c} 
\toprule
 & nref & drag max & lift max & time & nonlinear its&  linear its\\
\midrule
\multirow{4}{*}{RadauIIA(2)}&3&2.95&0.48&0.03&1.37&1.61 \\ 
                                &4&2.95&0.48&0.04&1.76&2.39 \\ 
                                &5&2.95&0.48&0.07&2.14&3.56 \\  
                                &6&2.95&0.48&0.13&2.52&5.09 \\ 
   
\midrule
\multirow{4}{*}{LobattoIIIC(2)} &3&2.95&0.48&0.03&1.35&1.80 \\ 
                                &4&2.95&0.48&0.05&1.79&2.76 \\ 
                                &5&2.95&0.48&0.07&2.15&4.30 \\  
                                &6&2.95&0.48&0.17&2.70&8.19 \\ 

\midrule
\midrule
\multirow{4}{*}{RadauIIA(3)}&3&2.95&0.48&0.04&1.57&2.27 \\ 
                            &4&2.95&0.48&0.07&1.90&2.78 \\ 
                            &5&2.95&0.48&0.11&2.17&3.68 \\  
                            &6&2.95&0.48&0.17&2.56&5.01 \\   

\midrule
\multirow{4}{*}{LobattoIIIC(3)}&3&2.95&0.48&0.04&1.38&1.74 \\ 
                               &4&2.95&0.48&0.06&1.78&2.50 \\ 
                               &5&2.95&0.48&0.10&2.14&3.60 \\  
                               &6&2.95&0.48&0.19&2.54&5.06 \\ 
\bottomrule
\end{tabular}

\caption{Maximum drag and lift values, average wall-clock time per time-step (in minutes) and average numbers of nonlinear and linear iterations per time step for $3\leq \ell \leq 6$ for Navier--Stokes flow past a cylinder.}
\label{2d NS results}
\end{table}
  
\begin{figure}
\centering
\includegraphics[width=0.45\textwidth]{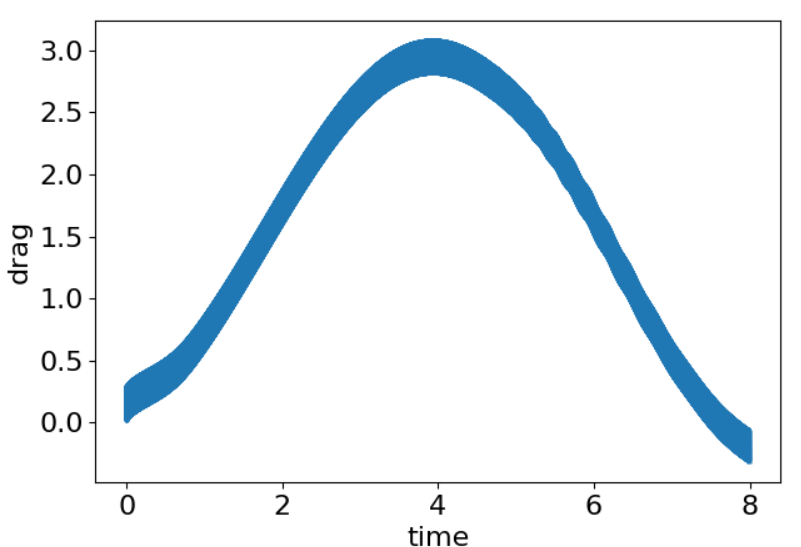}
%\subcaption{Drag}
%\label{fig:subfig1}
\qquad
\includegraphics[width=0.45\textwidth]{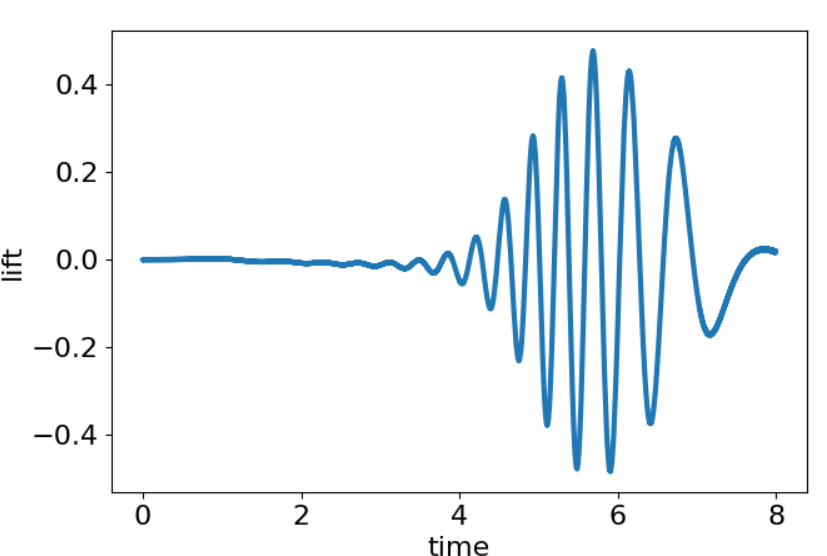}
%\subcaption{Lift}
%\label{fig:subfig2}
\caption{Drag and lift for $\ell=4$ using Gauss(3). The ``thick lines'' indicate that the solutions are highly oscillatory in time, due to the lack of L-stability of the integrator.}
\label{fig: NS GL3}
\end{figure}

Using both RadauIIA and LobattoIIIC IRK discretizations, with either two or three stages, results in computed lift and drag values that are consistent with those presented in \cite{john2004reference, Irksome}. However, results computed with Gauss were not. \Cref{fig: NS GL3} shows results using the three-stage Gauss method with $\ell=3$, computed with a stricter stopping tolerance (absolute nonlinear residual norm below $1/N^4$) than used above for RadauIIA and LobattoIIIC methods\footnote{Using the same stopping tolerance led to even more inconsistent data.}. The appearance of ``thick lines'' in these plots reflects highly oscillatory numerical solutions.  We hypothesize that this is due to the lack of L-stability of the integrator, where large negative eigenvalues of the linearized spatial operator are not quickly damped but, rather, slowly decay and oscillate in time due to a stability function value close to $-1$.  Refinement in time for fixed spatial grids should ameliorate the issue, but leads to increased computational costs to achieve similar accuracy to that given by RadauIIA and LobattoIIIC with these timesteps.

\subsection{Two-Dimensional MHD Island Coalescence}

We next consider a standard test model in MHD, of two-dimensional island coalescence.  This model mimics flow in a large aspect ratio tokamak, considering a cross-section of flow of magnetically confined plasma.  When a large external magnetic field is imposed in the ``toroidal'' direction of the tokamak, essentially two-dimensional dynamics result.  This model geometry is then mapped and rescaled to a square domain, $\Omega = (-1,1)^2$, with periodic boundary conditions on the left and right edges (see \cite{knoll2006coalescence, scottmhd, adler2013island} for more details). In this geometry, an equilibrium solution to the MHD equations is given by 
$$\begin{aligned}
\textbf{u}_{0}(x,y)&=\textbf{0},\\
\textbf{B}_{0} (x,y)& = \frac{1}{\cosh(2\pi y)+k\cos(2\pi x)}\begin{pmatrix} \sinh(2\pi y)  \\ k\sin(2\pi x)  \end{pmatrix} ,\\
p(x,y)&=\frac{1-k^{2}}{2}\left(1+\frac{1}{(\cosh(2\pi y)+k\cos(2\pi x))^{2}}\right),\\
\gamma(x,y)&=0,
\end{aligned}$$
 where $k=0.2$, when forcing terms of
$$\begin{aligned}
\textbf{f}_{\textbf{u}}& = \textbf{0},\\
\textbf{f}_{\textbf{B}}& = \frac{-8\pi ^{2}(k^{2}-1)}{Re_{m}( \cosh(2\pi y)+k\cos(2\pi x))^{3}}\begin{pmatrix} \sinh(2\pi y)  \\ k\sin(2\pi x)  \end{pmatrix},
 \end{aligned}$$
 are imposed on the differential equation.  To initialize a dynamic problem, these forcing terms are applied, but the initial condition is perturbed by adding 
$$\begin{aligned}
\delta \textbf{B}= \frac{-0.01}{\pi}\begin{pmatrix}
-\cos(\pi x)\sin(\frac{\pi y}{2})/2 \\
\cos(\frac{\pi y }{2})\sin(\pi x)
\end{pmatrix}.
\end{aligned}$$
 to the equilibrium solution at $t=0$.  The expected effect of this perturbation is to create two initially separated ``islands'' of current density that break the magnetic field lines, which then reconnect.
At the reconnection point (or $\mathcal{X}$-point), a sudden sharp spike should be seen in the magnetic current density.  At higher Reynolds numbers, a ``sloshing'' effect should occur before the islands of current density merge.

\begin{figure}
\centering
\includegraphics[width=0.95\textwidth]{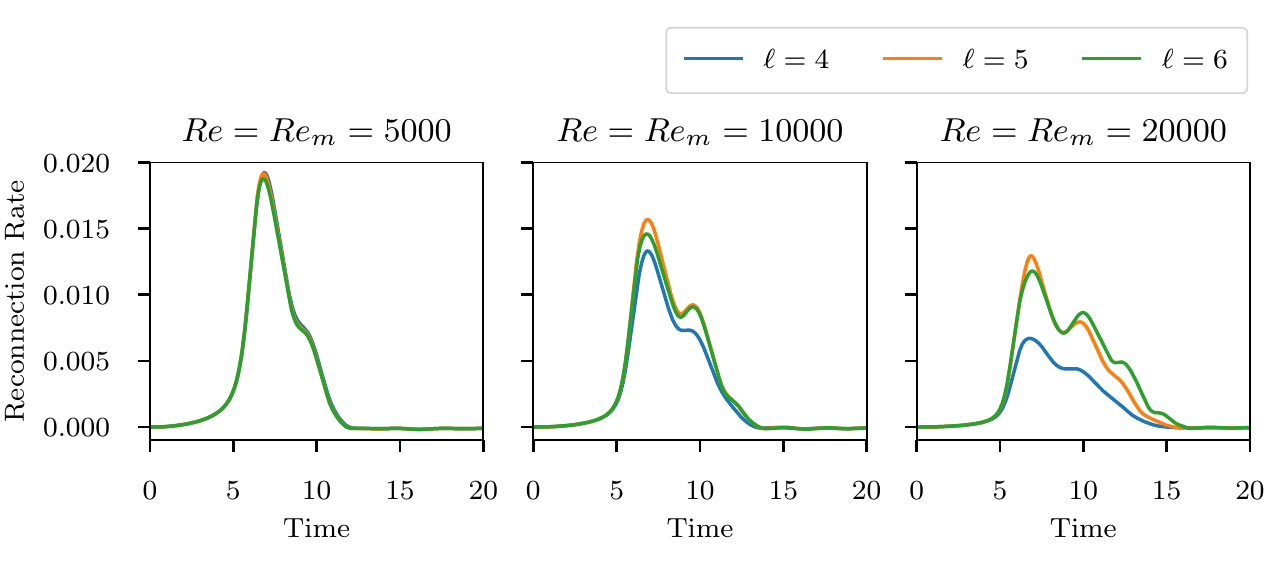}
\caption{Reconnection rates recorded  for the 2D MHD island coalescence model with varying Reynolds numbers of $\text{Re} = \text{Re}_{m} = 5000$, $10000$ and $20000$ on 3 different levels of refinement using the LobattoIIIC(2) temporal discretization.}
\label{fig:recon rate}
\end{figure}

As above, no analytical solution is known for this problem. A key measure of the physical fidelity is the time-history of the \emph{reconnection rate}, computed as the difference between the curl of $\textbf{B}$ at the origin at the current time and its value at the origin, scaled by $1/\sqrt{\text{Re}_m}$.  We compute this using the same methodology as in~\cite{scottmhd}.  As $\text{Re}$ and $\text{Re}_m$ increase, the peak value of the reconnection rate should decrease, and the length of time for which this value is nonzero should increase.  In this section, we consider only the two-stage LobattoIIIC integrator, and integrate until $T_f = 20$.  Following~\cite{scottmhd}, we ``substep'' for the first time-step, taking 10 substeps to initialize the simulation and avoid problems with nonlinear convergence. We consider a coarsest spatial mesh of $20\times 20$ quadrilateral elements, again each cut into 4 triangles, and present results for $\ell=4,5,6$ refinements.  For $\ell=4$, the resulting discretization has about 2.7 million DoFs per stage, while it has about 42.7 million DoFs per stage for $\ell=6$.  For $\ell=4$, we use $\Delta t = 0.025$, which is halved with each spatial refinement. We test for 3 different pairs of Reynolds numbers, $\text{Re}=\text{Re}_{m}$: $5000$, $1000$ and $20000$.  \Cref{fig:recon rate} shows the computed reconnection rates for these problems with varying $\text{Re} =\text{Re}_m$ and $\ell$, properly reflecting the expected behaviour.

\begin{figure}
\centering
\includegraphics[width=0.95\textwidth]{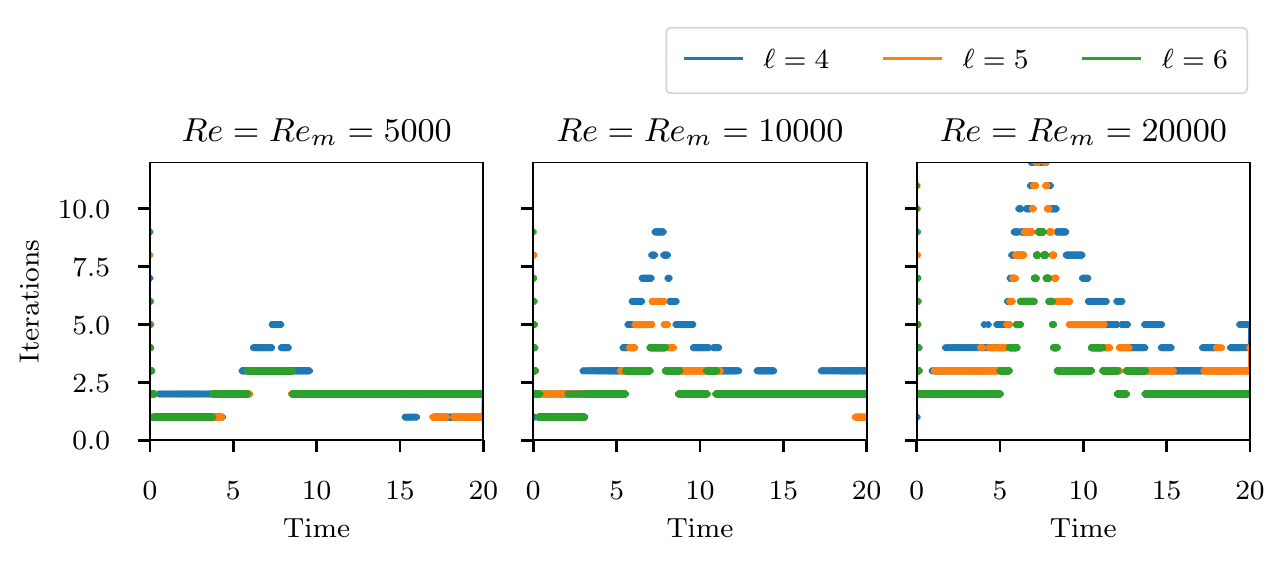}
\caption{Number of linear iterations per timestep for the 2D MHD island coalescence model with varying Reynolds numbers of $\text{Re} = \text{Re}_{m} = 5000$, $10000$ and $20000$ on 3 different levels of refinements using the LobattoIIC(2) integrator.}
\label{fig:linear its}
\end{figure}

\begin{figure}
\centering
\includegraphics[width=0.95\textwidth]{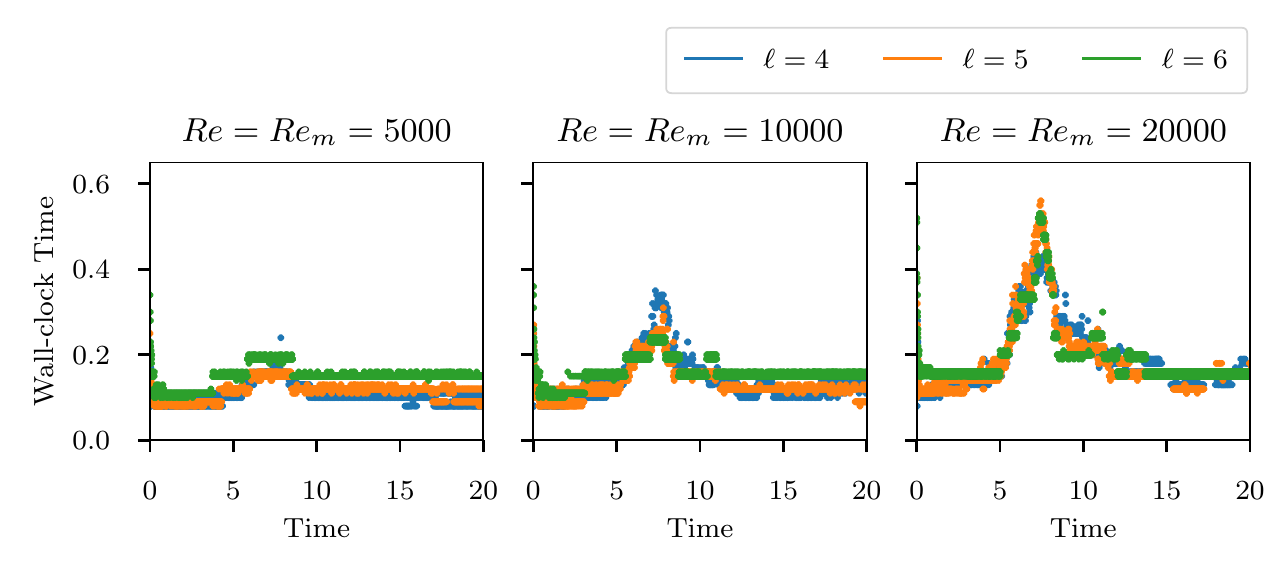}
\caption{Wall-clock time (in minutes) for the nonlinear system solve at each time-step for the 2D MHD island coalescence model with varying Reynolds numbers of $\text{Re} = \text{Re}_{m} = 5000$, $10000$ and $20000$ on 3 different levels of refinements using the LobattoIIC(2) integrator.}
\label{fig:nonlinear time}
\end{figure}

For this problem, we adjust the nonlinear and linear solver parameters as follows.  Taking $N=20\times 2^\ell$ as a representative number of elements in one dimension on refinement level $\ell$, we set both nonlinear and linear stopping tolerances to demand an absolute reduction of the $\ell_2$ norm of the corresponding residual below $1/N^2$.  We now use 3 pre- and post-relaxation sweeps, with the Chebyshev polynomials defining the relaxation taken over the interval [2, 10].  \Cref{fig:linear its,fig:nonlinear time} present results from a weak scaling study, using 40 cores on 1 node for $\ell=4$, 160 cores on 4 nodes for $\ell=5$, and 640 cores on 16 nodes for $\ell=6$.  We note that these are larger core counts than those used for the same underlying meshes with a BDF2 discretization in~\cite{scottmhd}; however, this is due to the larger number of DoFs in the system using a 2-stage discretization.  On average, our finest-grid problems have about 67 thousand DoFs per stage per core, which is a reasonable range for weak scaling.  We note that the $\ell=4$ problem takes about 2 hours of wall-clock time with these settings, with slightly better than doubling of wall-clock with each refinement (due to the halving of $\Delta t$ with each refinement, but also improved solver performance).

\Cref{fig:linear its} shows the number of linear solver iterations recorded per timestep as we vary $\text{Re} = \text{Re}_m$ and $\ell$.  We note slight growth in iteration counts as $\text{Re} = \text{Re}_m$ increases (and the problem becomes less diffusive in nature), but also improving iteration counts at fixed values of $\text{Re} = \text{Re}_m$ as $\ell$ increases.  Comparing with iteration counts from~\cite{scottmhd}, we see slightly higher iteration counts here, with slightly worse dependence on $\text{Re} = \text{Re}_m$, but still reasonable performance overall.  Wall-clock times per timestep, shown in \Cref{fig:nonlinear time}, generally reflect the linear iteration counts.  In particular, we again see a general increase with $\text{Re} = \text{Re}_m$, and a general decrease with increasing $\ell$.  The most expensive solves in the test set are still achieved in under 1 minute per time-step, and the average time is about 0.2 minutes per time-step.


\begin{figure}
\centering
\includegraphics[width=0.95\textwidth]{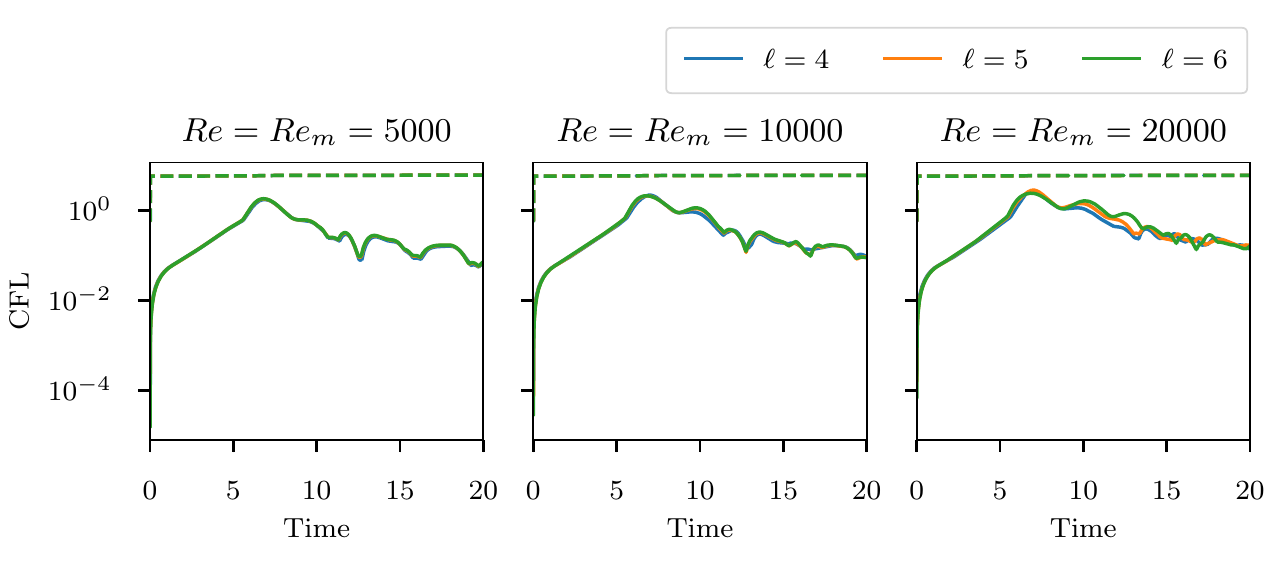}
\caption{CFL values at each timestep for the 2D MHD island coalescence model with varying Reynolds numbers of $\text{Re} = \text{Re}_{m} = 5000$, $10000$ and $20000$ on 3 different levels of refinements using the LobattoIIC(2) integrator. Solid lines represent fluid CFL, while dashed lines represent Alfv\'en CFL.}
\label{fig:cfl}
\end{figure}

To better understand the linear and nonlinear solver performance shown above, we compute both the fluid and magnetic (Alfv\'en) CFL numbers for the flow.  At each timestep, for the given solutions for $\textbf{u}$ and $\textbf{B}$, we approximate the maximum magnitude of the vector fields (by projecting $\textbf{u}\cdot\textbf{u}$ and $\textbf{B}\cdot\textbf{B}$ into the discontinuous piecewise-constant finite-element space on the finest mesh and computing the maximum values of these projections), $u_{\text{max}}$ and $B_{\text{max}}$, and then computing the fluid CFL value, $u_{\text{max}}\frac{\Delta t}{h}$, and the Alfv\'en CFL value, $B_{\text{max}}\frac{\Delta t}{h}$, where $h$ is a representative edge length for the spatial mesh. \Cref{fig:cfl} shows both CFL values calculated at each timestep for the simulations considered, showing identical results to those obtained in the BDF2 case in~\cite{scottmhd}.  We note that, aside from the initial substeps, the Alfv\'en CFL is roughly constant at a value around 6, while the fluid CFL peaks at the same time as the reconnection rate, and is above 1 for the largest values of $\text{Re} = \text{Re}_m$ considered.

\subsection{Three-dimensional MHD lid-driven cavity}

Finally, we present results for a three-dimensional lid-driven cavity MHD model on the unit cube, $\Omega = (0,1)^{3}$, following~\cite{phillips2016block}. On the top face, $z=1$, the flow is driven by imposed velocity $\textbf{u}=(1,0,0)^T$, while $\textbf{u}=(0,0,0)^T$ on all other faces. The tangential components of the magnetic field are set to match those of $\textbf{B}=(-1,0,0)^T$ on all faces of the cube. We set $\gamma=0$ on all boundary faces and fix the pressure $p=0$ at the origin.  In this section, we consider 3-grid methods, refining a given coarsest grid twice for each test.  This is driven by the consideration that, with increasing finest grid size, we require more cores over which to parallelize the computation; however, there is a software limitation within Firedrake that requires that the coarsest grid in the simulation must have at least 1 cell per core. While satisfying this requirement is not burdensome in 2D, it becomes problematic with increasing memory and computational requirements of 3D simulations. In all cases, we construct the coarsest grid by taking a uniform hexahedral mesh of the cube, then cutting each hexahedral element into 6 tetrahedra in the usual way.  We still use $\ell$ to denote the levels of refinement, but now $\ell=1$ denotes the smallest grid, created by refining a $2\times 2\times 2$ grid twice, while $\ell=2$ denotes the grid created by refining a $4\times 4\times 4$ grid twice, and $\ell=3$ denotes the grid created by refining a $8\times 8\times 8$ grid twice.  With $\ell=1$, our discretization has about 20 thousand DoFs per stage, increasing to about 1.1 million DoFs per stage for $\ell=3$.

We employ the same spatial discretization and again use the LobattoIIIC(2) integrator. We integrate until $T_f = 2.5$.  For $\ell=1$, we take $\Delta t = 0.125$, and halve $\Delta t$ with each refinement.  An all-zero initial condition is used.  \Cref{table: 3d LDC results} presents average linear and nonlinear iterations per timestep, along with average wall-clock time per nonlinear solve for the three grids above and $\text{Re} = \text{Re}_m = 10^p$ for $1 \leq p \leq 3$.  For $\ell=1$, 10 cores on 1 node are used, increasing to 80 cores on 2 nodes for $\ell=2$ and 640 cores on 16 nodes for $\ell=3$.  We use 3 pre- and post-relaxation sweeps, here accelerated using GMRES, as this was observed to result in better overall iteration counts and computation times than using Chebyshev acceleration, likely due to the convective nature of the problem at high Reynolds numbers. The nonlinear solve at each timestep requires the absolute value of the $\ell_2$ norm of the residual to be reduced below $10^{-6}$, and the same stopping criterion is used for the linear solves as well.

\begin{figure}
        \centering
            \includegraphics[width=0.475\textwidth]{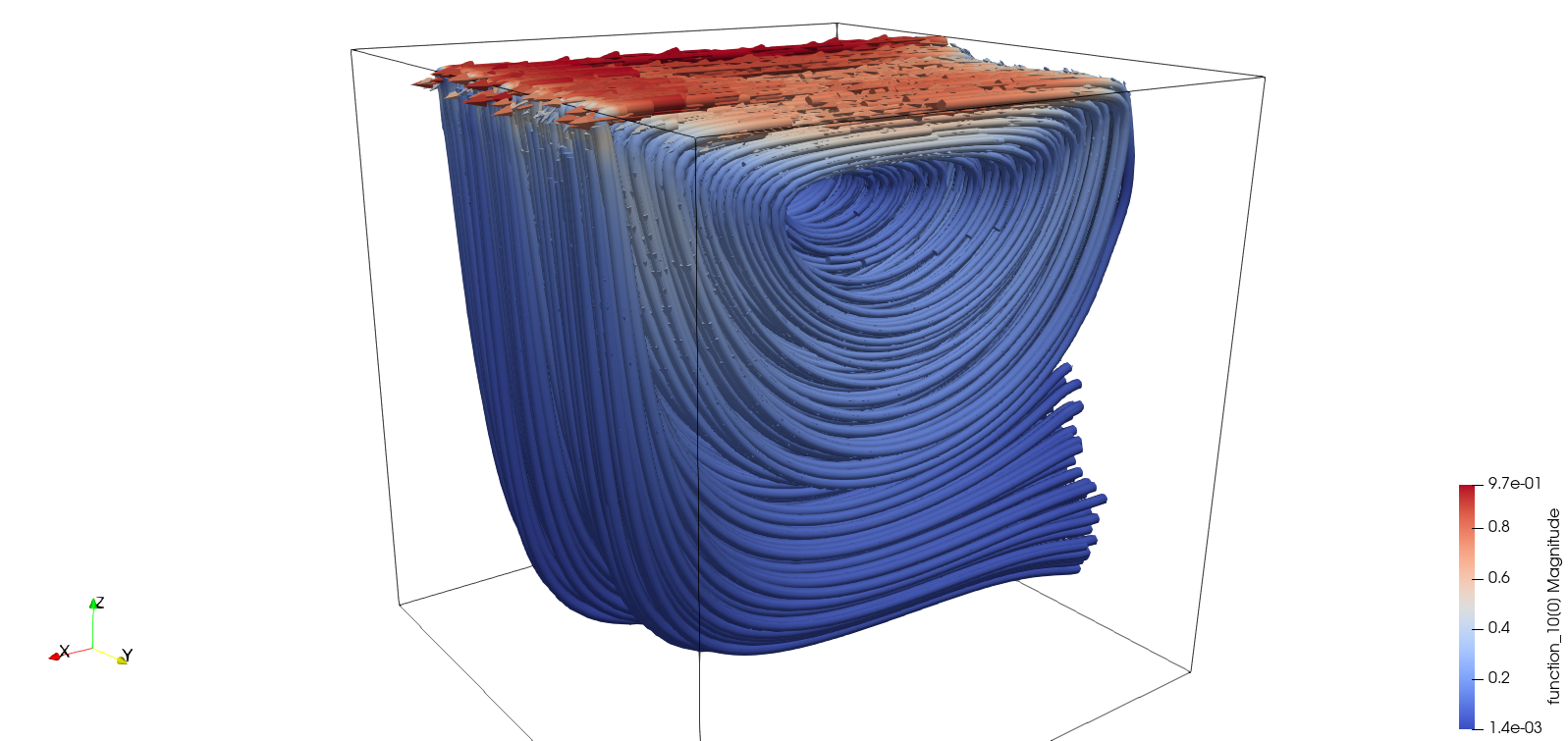}
            \includegraphics[width=0.475\textwidth]{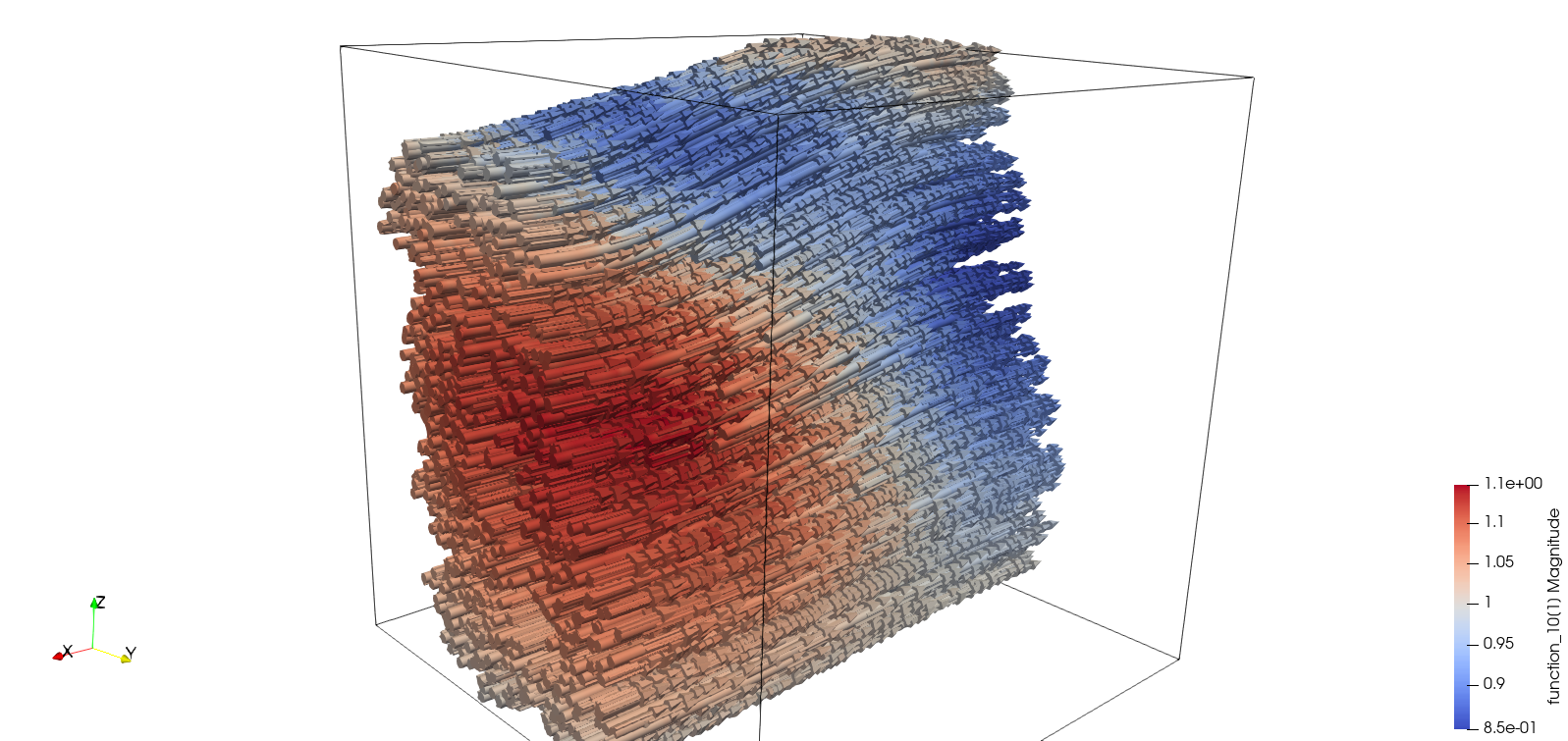}
        \caption{Streamlines of velocity (left) and magnetic field (right) for $\ell=3$ with $\text{Re} = \text{Re}_m = 10$.} 
        \label{fig: LDC}
      \end{figure}

Several trends can be observed in these results.  First, for fixed values of $\text{Re} = \text{Re}_m$, we generally observe improving solver performance as $\ell$ is increased, as expected.  Similarly, we typically observe degrading solver performance as $\text{Re} = \text{Re}_m$ is increased for fixed $\ell$.  Overall, the iteration counts are quite reasonable, except for $\text{Re} = \text{Re}_m = 1000$ with $\ell=1,2$.  Here, the problem is quite severely under-resolved, with a finest-grid mesh spacing of $h = 0.0625$ with $\ell=2$, so it is not surprising that the solver suffers when the discretization is so poor. For smaller Reynolds numbers, $\text{Re}=\text{Re}_{m}=1$ (not shown here), using Chebyshev acceleration gave significantly better results than using GMRES-accelerated relaxation, which failed to converge in some cases. \Cref{fig: LDC} presents representative solutions for $\ell=3$ with $\text{Re} = \text{Re}_m = 10$ (where the solutions are well-resolved), showing streamlines of both the velocity field, $\textbf{u}$, and the magnetic field, $\textbf{B}$, at the final time at refinement $\ell=3$.  

\begin{table}
\centering
\begin{tabular}{c|c|c|c|c} 
\toprule
$\ell$ &  & $\text{Re}=10$&$\text{Re}=100$&$\text{Re}=1000$\\
\midrule
\multirow{3}{*}{1} & linear its.     & 6.86 & 8.34 & 31.55 \\ 
                   & nonlinear its.  & 3.10 & 2.76 & 3.52 \\ 
                   & time            & 0.11 & 0.11 & 0.24  \\  
\midrule
\multirow{3}{*}{2} & linear its     &  6.61 & 5.41 & 17.59 \\ 
 				   & nonlinear its  &  3.18 & 2.43 & 3.02  \\
				   & time           &  0.21 & 0.17 & 0.30\\                                   
\midrule
\multirow{3}{*}{3} & linear its   & 5.76 & 4.20 & 6.08 \\
                   & nonlinear its& 2.57 & 2.33 & 2.16\\
                   & time         & 0.22 & 0.22 & 0.22\\   

\bottomrule
\end{tabular}
\caption{Average number of linear and nonlinear iterations per time-step and wall-clock time per nonlinear iteration (in minutes) for the 3D MHD lid-driven cavity problem with various Reynolds numbers and grid refinements, using the LobattoIIIC(2) integrator. }
\label{table: 3d LDC results}
\end{table}

\section{Conclusion}
\label{sec:conclusion}      

In this paper, we have developed monolithic Vanka relaxation schemes for fully-implicit Runge--Kutta discretizations of saddle point problems arising in models of fluid flow.
Within a Newton--Krylov--multigrid setting, our method is shown to be effective for both Newtonian and magnetohydrodynamic flows, in both two and three spatial dimensions. The algorithm is chosen with parallel implementation in mind, and weak scaling results are shown up to 640 cores.

There are many possibilities for future work.  We note primarily that the current study uses relatively low-order spatial discretizations, based on classical Taylor--Hood elements for velocity and pressure. A next step in this research is to extend these solvers to more sophisticated finite-element discretizations that preserve the incompressibility and solenoidality constraints exactly, as in~\cite{laakmann2021augmented, KHu_etal_2017a}.  An important question for future work is the extension of these techniques to higher-order discretizations, where the cost of classical sparse direct solvers for the patch problems becomes prohibitive.

%\section*{References}

\bibliography{references}

\end{document}